\documentclass[11pt]{article}
\usepackage{amssymb}
\usepackage{graphicx}
\usepackage{color}

\usepackage{amsmath}

\usepackage{amsfonts}

\usepackage{amssymb}
\usepackage{graphics}
 \usepackage{bbm}

\usepackage{titlesec}
\usepackage[titletoc,toc,title]{appendix}

\titleformat{\section}{\large\bfseries}{\thesection.}{0.3em}{}
\titleformat{\subsection}{\bfseries}{\thesubsection.}{0.3em}{}
\titleformat{\subsubsection}{\small\bfseries}{\thesubsubsection.}{0.3em}{}

 \allowdisplaybreaks[2]

\usepackage{latexsym}

\renewcommand{\skew}{\mathop{\rm skew}}
 \DeclareMathOperator{\sym}{sym}
\DeclareMathOperator{\tr}{tr} 
 \DeclareMathOperator{\dev}{dev}
\DeclareMathOperator{\Curl}{Curl\hskip.04truecm}
\DeclareMathOperator{\Curlt}{Curl\,Curl}
\DeclareMathOperator{\Div}{Div} \DeclareMathOperator{\curl}{curl\,}

\DeclareMathOperator{\sL}{\mathfrak{sl}}
\DeclareMathOperator{\so}{\mathfrak{so}}

\newcommand{\Chi}{\raisebox{0.5ex}{\mbox{{\Large $\chi$}}}}

\newcommand{\yieldlimit}{{\sigma}_{\mathrm{y}}}

\newcommand{\yieldzero}{{\sigma}_{0}}

\newcommand{\C}{\mathbb{C}}

\newcommand{\PL}{\mathcal{X}}
\newcommand{\rgam}{\raisebox{0.2ex}{$\gamma$}}
\newcommand{\reta}{\raisebox{0.2ex}{$\eta$}}
\newcommand{\vectgam}{\raisebox{0.2ex}{\underline{$\gamma$}}}
\newcommand{\vecteta}{\raisebox{0.2ex}{$\underline{\eta}$}}
\newcommand{\QL}{\mbox{\footnotesize Q}}
\newcommand{\BBR}{\mbox{$\mathbb{R}$}}

\newcommand{\SFH}{\mbox{$\mathsf{H}$}}

\newcommand{\SFQ}{\mbox{$\mathsf{Q}$}}
\newcommand{\SFV}{\mbox{$\mathsf{V}$}}
\newcommand{\SFW}{\mbox{$\mathsf{W}$}}
\newcommand{\SFP}{\mbox{$\mathsf{P}$}}
\newcommand{\SFZ}{\mbox{$\mathsf{Z}$}}

\newcommand{\bvarepsilon}{\mbox{$\bf\varepsilon$}}

\newcommand{\dsize}{\displaystyle}
\renewcommand{\div}{\mathop{\rm div}\nolimits}
\newcommand{\ba}{\mbox{\boldmath{$a$}}}
\numberwithin{equation}{section}
\newcommand{\parent}[3]{\left #1 {#3} \right #2}
\newcommand{\graffe}[1]{\parent \{ \}{#1}}

\newcommand{\norm}[1]{\lVert{#1}\rVert} 
\newcommand{\bfig}[2]{\begin{figure}\begin{center}\begin{picture}(341.8,#2)(
#1,0)}
\newcommand{\efig}[2]{\end{picture}\caption{#2.}\lbl{#1}\end{center}
\end{figure}}

\newcommand{\la}{\langle}
\newcommand{\ra}{\rangle}
\newcommand{\id}{{\boldsymbol{\mathbbm{1}}}}
\newcommand{\qed}{\qquad$\blacksquare$}

\newtheorem{lemma}{Lemma}[section]

\newtheorem{remark}{Remark}[section]

\newcommand{\be}{\begin{equation}}

\newcommand{\ee}{\end{equation}}


 \makeatletter
 \let\@fnsymbol\@arabic
 \makeatother
\begin{document}
\vskip-3truecm 
\title{
\vspace{-1.in} {\Large Well-posedness for the microcurl model in both single and  polycrystal gradient plasticity
\texttt{}}}

\author{{\large Fran\c{c}ois Ebobisse{\footnote{Corresponding author, Fran\c{c}ois Ebobisse,
Department of Mathematics and Applied Mathematics, University of
Cape Town, Rondebosch 7700, South Africa, e-mail:
francois.ebobissebille@uct.ac.za}}\quad and\quad Patrizio Neff{\footnote
{Patrizio Neff, Lehrstuhl f\"ur  Nichtlineare Analysis und
Modellierung, Fakult\"{a}t f\"ur Mathematik, Universit\"at
Duisburg-Essen, Thea-Leymann Str. 9, 45127 Essen, Germany, e-mail:
patrizio.neff@uni-due.de, http://www.uni-due.de/mathematik/ag\b{
}neff}} \quad and\quad  Samuel Forest{\footnote {Samuel Forest, MINES ParisTech, Centre des Mat\'eriaux, UMR CNRS 7633, BP 87, 91003 Evry, France, e-mail: samuel.forest@ensmp.fr}} }\vspace{1mm}}



\date{\today}

 \maketitle

\begin{center}

 {{\small
\begin{abstract}
 We consider the recently introduced microcurl model which is a variant of strain gradient plasticity in which the curl of the plastic distortion is coupled to an additional micromorphic-type field. For both single crystal and polycrystal cases, we formulate the model and show its well-posedness in the rate-independent case provided some  local hardening (isotropic or linear kinematic) is taken into account. To this end, we use the functional analytical framework developed by Han-Reddy. We also compare the model to the relaxed micromorphic model as well as to a dislocation-based gradient plasticity model.

   \end{abstract}}}
\end{center}
\noindent {\bf Key words:} plasticity, gradient plasticity, variational modeling, dissipation function, micromorphic continuum, defect energy, micro-dislocation.
 
\vskip.2truecm\noindent {\bf AMS 2010 subject classification:}
35D30, 35D35, 74C05, 74C15, 74D10, 35J25.
\newpage
 {\small \tableofcontents}
\section{Introduction}\label{Intro}
In this paper we consider the so-called microcurl model in plasticity. The model was introduced in  \cite{Forest2013,CORFORBUS} to serve the purpose of augmenting classical plasticity with length scale effects while otherwise keeping the algorithmic structure of classical plasticity. The idea is simple and straight-forward: the in general non-symmetric plastic distortion $p$ (single crystal plasticity and polycrystal plasticity with plastic
spin) with its local in space evolution is energetically coupled to a micromorphic-type additional non-symmetric tensorial variable $\PL_p$ via a penalty-like term $\frac12\,\mu\,H_\chi\,\norm{p-\PL_p}^2$. The tensor field $\PL_p$ is generally assumed to be incompatible i.e.,  it may not derive from a vector field. The total energy in the model is then augmented by a quadratic contribution acting on the Curl of $\PL_p$. The new variable $\PL_p$ is now determined by free-variation of the energy w.r.t. the displacement field $u$ and the micromorphic field $\PL_p$ together with corresponding tangential boundary conditions for $\PL_p$. This generates the usual equilibrium equation on the one hand and  what we will call a micro-balance equation for $\PL_p$ on the other hand.

In the penalty limit $H_\chi\to\infty$, when one expects $p=\PL_p$, the variable $\Curl\PL_p$ is then interpreted to be the dislocation density tensor $\Curl p$. The advantage of such a formulation is clear: there is no need for an extended thermodynamic setting, since $\PL_p$ is not directly taking part in the dissipation. Constitutive laws including dissipative contributions of the microdeformation and microcurl can be proposed, as done
in \cite{FORSIEV2006} in the general micromorphic case, but they will require additional material parameters whose identification necessitates
material specific physical considerations. Thus, also no higher order boundary conditions at interfaces between elastic and plastic parts need be discussed. The resulting model can therefore be described as a {\it pseudo-regularized strain gradient model}.

The microdeformation variable $\PL_p$ has at least two different interpretations. 
First, it can be regarded as a mathematical auxiliary variable used to replace the higher order
partial differential equations arising in strain gradient plasticity by a system of two sets of
second order partial differential equations for the displacement and microdeformation. This method has 
computational merits for the implementation of strain gradient plasticity models in finite element codes,
see \cite{anandoz12}. In that case, $H_\chi$ is regarded as a mere penalty parameter and should be large enough to enforce the constraint
$p = \PL_p$. In contrast, the microdeformation variable $\PL_p$ can also be viewed as a
constitutive variable with physical interpretation, for instance based on statistical mechanics,
$p$ representing the average plastic distortion over the material volume element and $\PL_p$ being related to the 
variance of plastic deformation inside this volume. This interpretation is similar to the microconcentration variable
introduced in  \cite{ubachs04,cism11joerg} to solve Cahn--Hilliard equations.
In this context, $H_\chi$ must be regarded as a true material higher order modulus to be identified from 
suitable experimental results.
Compared to standard strain gradient plasticity, the microcurl model therefore possesses one additional parameter,
$H_\chi$, which allows for better description of physical results, as suggested in \cite{CGFBGK2010}.
An interpretation of the microdeformation $\PL_p$ was recently proposed in the case of polycrystalline plasticity and damage in
\cite{poh13homog,poh16} where it is related to the grain to grain heterogeneitiy of plastic deformation. \\
Another computational advantage of the microcurl single crystal plasticity model is that the number of independent degrees of freedom
($9$ tensor components of $\PL_p$, or $8$ in the case of incompressible microplasticity) is independent of the crystallographic
structure of the material and of the number of slip systems. This is in contrast to strain gradient plasticity models
involving the directional gradient of the slip variables \cite{GURT2002}, which require as many degrees of freedom as slip systems ($12$ at least
in FCC crystals, up to $48$ in BCC crystals!). A comparison and discussion of models based on the full dislocation density tensors with 
models involving densities of geometrically necessary dislocations can be found in \cite{mesarovic15}.

In this paper we will consider the microcurl model in two variants. First, in its original form as a computational approach towards single crystal strain gradient plasticity. We formulate the governing system and show its well-posedness in the rate-independent case. The natural solution space for the micro-variable $\PL_p$ is  the Sobolev-like space $H(\Curl)$.

Second, we extend the approach formally to polycrystalline plasticity in which the plastic variable $\bvarepsilon_p:=\sym p$ (the plastic strain) is assumed to be symmetric. In this case we still allow for a non-symmetric micro-variable $\PL_p$ which is now coupled to the plastic variable only via its symmetric part by $\frac12\,\mu\,H_\chi\,\norm{\sym(p-\PL_p)}^2$. This represents an alternative to recently proposed strain gradient 
plasticity models involving a plastic spin tensor for polycrystals \cite{GURT2004,bardella15torsion,poh16}. Again, we show the well-posedness of the formulation. Here, we need recently introduced coercive inequalities generalizing Korn's inequality to incompatible tensor fields \cite{NPW2011-1, NPW2012-1, NPW2012-2,NPW2014}. 

The mathematical analysis (with results of existence and uniqueness) for both variants (single crystal and polycrystalline) is obtained through the machinery developped by Han-Reddy \cite[Theorems 6.15 and 6.19]{Han-ReddyBook} for classical plasticity and recently extended to models of gradient plasticity in \cite{DEMR1, NCA, EBONEFF, ENR2015, EHN2016}. In this  approach, the model through the primal form of the flow rule is weakly formulated as a variational inequality and the key issue for its well-posedness is the study of the coercivity of the  bilinear form involved on a suitable closed convex subset of some Hilbert space.

The polycrystalline variant of the microcurl model bears some superficial resemblence with the recently introduced relaxed micromorphic models \cite{NGMPR2014, NGMP2014}. For purpose of clarification, we present the relaxed micromorphic model and clearly point out the differences. In order to put the microcurl modelling framework further on display we finish this introduction with another dislocation based strain gradient plasticity model with plastic spin
 \cite{EBONEFF, ENR2015, EHN2016}, see Table \ref{table:kinhardspin}. 

\begin{table}[ht!]{\footnotesize\begin{center}
\begin{tabular}{|ll|}\hline &\qquad\\
 {\em Additive split of distortion:}& $\nabla u =e +p$,\quad $\bvarepsilon^e:=\sym e$,\quad $\bvarepsilon^p:=\sym p$\\
{\em Equilibrium:} & $\mbox{Div}\,\sigma +f=0$ with
$\sigma=\C_{\mbox{\tiny iso}}\bvarepsilon^e$\\ {\em Free energy:} &
$\frac12\,\langle\C_{\mbox{\tiny iso}}\bvarepsilon^e,\bvarepsilon^e\rangle+\frac12\,\mu\, k_1\,\norm{\bvarepsilon^p}^2+\frac12\,\mu\,
L^2_c\,\norm{\Curl p}^2$\\&\\
{\em Yield condition:} &
$\phi(\Sigma_E):=\norm{\dev\Sigma_E}-\yieldzero\leq0$\\
 {\em where } & $\Sigma_E:=\sigma+\Sigma^{\mbox{\scriptsize
lin}}_{\mbox{\scriptsize curl}}+\Sigma^{\mbox{\scriptsize
lin}}_{\mbox{\scriptsize kin}}$,\,\, \\ &$\Sigma^{\mbox{\scriptsize
lin}}_{\mbox{\scriptsize curl}}:=-\mu\, L^2_c\,\Curlt p$,\,\quad
$\Sigma^{\mbox{\scriptsize lin}}_{\mbox{\scriptsize kin}}:=-\mu\,
k_1\,\bvarepsilon^p$
 \\&\\{\em Dissipation inequality:} &
 $\dsize\int_\Omega\langle\Sigma_E,\dot{p}\rangle dx\geq0$\\
 {\em Dissipation function:} &$\mathcal{D}(q):=\yieldzero \norm{q}$\\
 {\em Flow rule in primal form:} &
 $\Sigma_E\in\partial \mathcal{D}(\dot{p})$\\&\\
{\em Flow rule in dual form:}
&$\dot{p}=\lambda\,\dsize\frac{\dev\Sigma_E}{\norm{\dev\Sigma_E}},\quad\qquad \lambda=\norm{\dot{p}}$\\&\\
{\em KKT conditions:} &$\lambda\geq0$, \quad $\phi(\Sigma_E)\leq0$,
\quad $\lambda\,\phi(\Sigma_E)=0$\\
 {\em Boundary conditions for $p$:} & $p\times{n}=0$ on
 $\Gamma_{\mbox{\scriptsize D}}$,\,\, $(\Curl p)\times{n}=0$ on $\partial\Omega\setminus\Gamma_{\mbox{\scriptsize D}}$\\
 {\em Function space for $p$:} & $p(t,\cdot)\in \mbox{H}(\mbox{Curl};\;\Omega,\,\BBR^{3\times 3})$\\
 \hline
\end{tabular}\caption{\footnotesize The polycrystalline plasticity model model with linear kinematic hardening and plastic
spin studied in \cite{EHN2016}.}\label{table:kinhardspin}\end{center}}\end{table}
Here, the microcurl-type regularization would be obtained by considering the microcurl energy 
\begin{equation}\label{micropoly-spin}\frac12\,\langle\C_{\mbox{\scriptsize iso}}\,\bvarepsilon^e,\bvarepsilon^e\rangle+\frac12\,\mu\,H_\chi\,\norm{p-\PL_p}^2+\frac12\mu\,k_1\norm{\sym p}^2+\frac12\,\mu\,
L^2_c\,\norm{\Curl \PL_p}^2 \end{equation} and for  $H_\chi\to\infty$ we would recover the model from Table \ref{table:kinhardspin}. 

The polycrystalline microcurl variant which we introduce in this paper is, however, based on the energy 
\begin{equation}\label{variant-micropoly-spin}\frac12\,\langle\C_{\mbox{\scriptsize iso}}\,\bvarepsilon^e,\bvarepsilon^e\rangle+\frac12\,\mu\,H_\chi\,\norm{\sym(p-\PL_p)}^2+\frac12\mu\,k_1\norm{\sym p}^2+\frac12\,\mu\,
L^2_c\,\norm{\Curl \PL_p}^2\,.\end{equation} 
This ansatz seems to be appropriate for polycrystalline plasticity without plastic spin.
\section{Some notational agreements and
definitions}\label{Notations} Let $\Omega$ be a bounded domain
 in $\BBR^3$ with Lipschitz continuous boundary $\partial\Omega$, which is occupied by the elastoplastic
body in its undeformed configuration. Let $\Gamma_{\mbox{\scriptsize D}}$ be a smooth
subset of $\partial\Omega$ with non-vanishing $2$-dimensional
Hausdorff measure. A material point in $\Omega$ is denoted by $x$
and the time domain under consideration is the interval $[0,T]$.\\
 For every $a,\,b\in\BBR^3$, we let $\la a,b\ra_{\BBR^3}$ denote the scalar
 product on $\BBR^3$ with associated vector
norm $|a|^2_{\BBR^3} = \la a, a\ra_{\BBR^3}$. We denote by
$\BBR^{3\times 3}$ the set of real $3\times 3$ tensors. The standard
Euclidean scalar product on $\BBR^{3\times 3}$ is given by $\la
A,\,B\ra_{\BBR^{3\times 3}} = \mbox{tr}\,\bigl[AB^T\bigr]$, where
$B^T$ denotes the transpose tensor of $B$. Thus, the Frobenius
tensor norm is $\norm{A}^2 = \la A,\,A\ra_{\BBR^{3\times 3}}$.
 In the following we omit the subscripts $\BBR^3$ and $\BBR^{3\times 3}$. The identity tensor on $\BBR^{3\times 3}$ will be denoted by
  $\id$, so that $\mbox{tr}(A) = \la A, \id\ra$.  The set
$\so(3):=\{X\in\BBR^{3\times 3}\,|\,\,X^T=-X\}$ is the Lie-Algebra
of skew-symmetric tensors.
 We let
$\mbox{Sym\,}(3):=\{X\in\BBR^{3\times 3}\,|\,\,X^T=X\}$ denote the
 vector space of symmetric tensors and $\sL(3):=\{X\in\BBR^{3\times
3}\,|\,\,\mbox{tr\,}(X)=0\}$ be the Lie-Algebra of traceless
tensors. For every $X\in\BBR^{3\times 3}$, we set
$\sym(X)=\frac12\bigl(X+X^T\bigr)$,
$\skew\,(X)=\frac12\bigl(X-X^T\bigr)$ and
$\dev(X)=X-\frac13\mbox{tr}\,(X)\,\id\in\sL(3)\,$ for the symmetric
part, the skew-symmetric part and the deviatoric part of $X$,
respectively. Quantities which are constant in space will be denoted
with an overbar, e.g., $\overline{A}\in\so(3)$ for the function
$A:\mathbb{R}^3\to\so(3)$ which is constant with constant value
$\overline{A}$.

The body is assumed to undergo infinitesimal deformations. Its
behaviour is governed by a set of equations and constitutive
relations. Below is a list of variables and parameters used
throughout the paper with their significations:\begin{itemize}
\item[$\bullet$] $u$  is the displacement of the macroscopic material
points;

\item[$\bullet$] $p$ is the infinitesimal plastic distortion variable which is a
non-symmetric second order tensor, incapable of sustaining
volumetric changes; that is, $p\in\sL(3)$. The tensor $p$\,
represents the average plastic slip; $p$ is not a state-variable, while the rate $\dot{p}$ is an infinitesimal state variable  in some suitable sense;

\item[$\bullet$] $e=\nabla u -p$ is  the infinitesimal elastic distortion which  is in general a
non-symmetric second order tensor and is an infinitesimal state-variable;

\item[$\bullet$] $\bvarepsilon_p=\sym p$ is the symmetric infinitesimal plastic strain
tensor, which is  trace free, {$\bvarepsilon_p\in\sL(3)$;} $\bvarepsilon_p$ is not a state-variable;
the rate $\dot{\bvarepsilon}_p$ is an infinitesimal  state-variable;

\item[$\bullet$] $\bvarepsilon_e=\sym\nabla u -\bvarepsilon_p$ is the symmetric infinitesimal  elastic
strain tensor and is an infinitesimal  state-variable;

\item[$\bullet$] $\PL_p\in\mathbb{R}^{3\times 3}$ is the non-symmetric 
infinitesimal  micro-distortion with $\sym\PL_p$ being the symmetric micro-strain;

\item[$\bullet$] $\sigma$   is the Cauchy stress tensor which is a symmetric
second order tensor and is an infinitesimal  state-variable;

\item[$\bullet$] $\yieldzero$ is the initial
yield stress for plastic variables $p$ or $\bvarepsilon_p:=\sym p$ and is an infinitesimal  state-variable;

\item[$\bullet$] $f$ is the body force;

\item[$\bullet$] $\Curl p=\alpha$ is the dislocation density
tensor satisfying the so-called Bianchi identities $\Div\alpha=0$ and is an infinitesimal  state-variable;

\item[$\bullet$] $\eta_p=\dsize\int_0^t\lVert\dot{\bvarepsilon}_p\rVert\,ds$ is
the accumulated equivalent plastic
strain and is an infinitesimal state-variable;

\item[$\bullet$] $\gamma^\alpha$ is the slip in the $\alpha$-th slip system in single crystal plasticity while $l^\alpha$ is the slip direction and $\nu^\alpha$ is the normal vector to the slip plane with $\alpha=1,\ldots,n_{\mbox{\tiny slip}}$. Hence, $p=\dsize\sum_\alpha\gamma^\alpha\,l^\alpha\otimes\nu^\alpha$.
\end{itemize}
\vskip.2truecm\noindent For isotropic media, the fourth order
isotropic elasticity tensor $\C_{\mbox{\scriptsize
iso}}:\mbox{Sym}(3)\to\mbox{Sym}(3)$ is given by
\begin{equation}
\C_{\mbox{\scriptsize iso}}\sym X = 2\mu\,\dev\,\sym X+\kappa
\,\tr(X) \id =2\mu\,\sym X+\lambda\,\tr(X)\id\label{C}
\end{equation}
for any second-order tensor $X$, where $\mu$ and $\lambda$ are the
Lam{\'e} moduli satisfying
\begin{equation}\label{Lame-moduli}
\mu>0\quad\mbox{ and }\quad 3\lambda +2\mu>0\,,
\end{equation} and $\kappa>0$ is the bulk modulus.
These conditions suffice for pointwise positive definiteness of the
elasticity tensor in the sense that there exists a constant $m_0 >
0$ such that
\begin{equation}
\forall X\in\mathbb{R}^{3\times 3}\mbox{ :}\,\quad\la \sym X,\C_{\mbox{\scriptsize iso}}\sym X\ra \geq m_0\,
\lVert\sym X\rVert^2\,. \label{ellipticityC}
\end{equation}

\vskip.2truecm\noindent The space of square integrable functions is
$L^2(\Omega)$, while the Sobolev spaces used in this paper are:
\begin{eqnarray}\label{sobolev-spaces}
\nonumber \mbox{H}^1(\Omega)&=&\{u\in
L^2(\Omega)\,\,|\,\,\mbox{grad\,}u\in
L^2(\Omega)\}\,,\qquad\quad\mbox{ grad}\,=\nabla\,,\\
\nonumber
&&\norm{u}^2_{H^1(\Omega)}=\norm{u}^2_{L^2(\Omega)}+\norm{\mbox{grad\,}u}^2_{L^2(\Omega)}\,,\qquad\forall u\in\mbox{H}^1(\Omega)\,,\\
 \mbox{H}(\mbox{curl};\Omega)&=&\{v\in
L^2(\Omega)\,\,|\,\,\mbox{curl\,}v\in
L^2(\Omega)\}\,,\qquad\mbox{curl\,}=\nabla\times\,,\\
\nonumber &&\norm{v}^2_{\mbox{\scriptsize
H}(\mbox{curl};\Omega)}=\norm{v}^2_{L^2(\Omega)}+\norm{\mbox{curl\,}v}^2_{L^2(\Omega)}\,,\,\,\quad\forall
v\in\mbox{H}(\mbox{curl;\,}\Omega)\,.
\end{eqnarray}
For every $X\in C^1(\Omega,\,\BBR^{3\times 3})$ with rows
$X^1,\,X^2,\,X^3$, we use in this paper the definition of $\Curl X$
in \cite{NCA, SVEN}:
\begin{equation}\label{def-Curl}\Curl X =\left(\begin{array}{l}\mbox{curl\,}X^1\,\,-\,\,-\\
\mbox{curl\,}X^2\,\,-\,\,-\\
\mbox{curl\,}X^3\,\,-\,\,-\end{array}\right)\in\BBR^{3\times
3}\,,\end{equation} for which $\Curl\,\nabla v=0$ for every $v\in
C^2(\Omega,\,\BBR^3)$. Notice that the definition of $\Curl X$
above is such that $(\Curl X)^Ta=\mbox{curl\,}(X^Ta)$ for every
$a\in\BBR^3$ and this clearly corresponds to the transpose of the
Curl of a tensor as defined in
\cite{GURTAN2005, GURTAN-BOOK}.\\

The following function spaces and norms will also be used later.
\begin{eqnarray}\label{Curl-spaces}
\nonumber \mbox{H}(\mbox{Curl};\,\Omega,\,\BBR^{3\times
3})&=&\Bigl\{X\in L^2(\Omega,\,\BBR^{3\times
3})\,\,\bigl|\,\,\mbox{Curl\,}X\in
L^2(\Omega,\,\BBR^{3\times 3})\Bigr\}\,,\\
\norm{X}^2_{\mbox{\scriptsize H}(\mbox{\scriptsize
Curl};\Omega)}&=&\norm{X}^2_{L^2(\Omega)}+\norm{\mbox{Curl\,}X}^2_{L^2(\Omega)}\,,\quad\forall
X\in\mbox{H}(\mbox{Curl;\,}\Omega,\,\BBR^{3\times 3})\,,\\
\nonumber
\mbox{H}(\mbox{Curl};\,\Omega,\,\mathbb{E})&=&\Bigl\{X:\Omega\to\mathbb{E}\,\,\bigl|\,\,X\in
\mbox{H}(\mbox{Curl};\,\Omega,\,\BBR^{3\times 3})\Bigr\}\,,
\end{eqnarray}
for $\mathbb{E}:=\sL(3)$ or
$\mbox{Sym}\,(3)\cap\sL(3)$.\vskip.2truecm\noindent
 We also consider the space
\begin{equation}\label{spacep-bc}\mbox{H}_0(\mbox{Curl};\,\Omega,\,\Gamma_{\mbox{\scriptsize D}},\BBR^{3\times 3})\end{equation} as
the completion in the norm in  (\ref{Curl-spaces}) of the space
$\bigl\{X\in C^\infty(\Omega,\,\BBR^{3\times
3})\,\,\bigl|\,\,\,X\times\,n|_{\Gamma_{\mbox{\scriptsize D}}}=0\bigr\}\,.$ Therefore, this
space generalizes the tangential Dirichlet boundary condition
$$X\times\,n|_{\Gamma_{\mbox{\scriptsize D}}}=0\,$$
to be satisfied by the plastic micro-distortion $\PL_p$. The space
$$\mbox{H}_0(\mbox{Curl};\,\Omega,\,\Gamma_{\mbox{\scriptsize D}},\mathbb{E})$$ is defined
as
 in (\ref{Curl-spaces}). \vskip.2truecm\noindent
 The divergence operator Div on second order
tensor-valued functions is also defined row-wise as
\begin{equation}\label{def-div}\mbox{Div}\,X=\left(\begin{array}{l}\mbox{div\,}X_1\\
\mbox{div\,}X_2\\
\mbox{div\,}X_3\end{array}\right)\,.\end{equation}

\section{The microcurl model in single crystal  gradient plasticity}\label{single}

\subsection{Kinematics}
Single-crystal plasticity is based on the assumption that the plastic deformation happens through crystallograpic shearing which represents the dislocation motion along  specific slip systems,  each being characterized
by a plane with unit normal $\nu^\alpha$ and slip direction $l^\alpha$ on that plane, and slips $\gamma^\alpha$ ($\alpha = 1, \ldots ,n_{\mbox{\tiny slip}}$). The flow rule
for the plastic distortion $p$ is written at the slip system level by means of the orientation tensor $m^\alpha$ defined as
\begin{equation}\label{orient-tensor} m^\alpha:=l^\alpha\otimes \nu^\alpha\,.\end{equation}
Under these conditions the plastic distortion $p$  takes the form
\begin{equation}\label{orientens} 
p=\sum_{\alpha=1}^{n_{\mbox{\tiny slip}}}\gamma^\alpha \,m^\alpha\end{equation}
so that the plastic strain $\bvarepsilon_p=\sym p$ is
\begin{equation}\label{orientens2}\bvarepsilon_p=\sum_{\alpha=1}^{n_{\mbox{\tiny slip}}} \gamma^\alpha \sym(m^\alpha)=\frac12\sum_{\alpha=1}^{n_{\mbox{\tiny slip}}}      \gamma^\alpha(l^\alpha\otimes \nu^\alpha+\nu^\alpha\otimes l^\alpha)\end{equation} and $\tr(p)=0$ since $l^\alpha\perp \nu^\alpha$.\vskip.2truecm\noindent
For the slips $\gamma^\alpha$ ($\alpha=1,\dots,n_{\mbox{\tiny slip}}$) we set
$$\mbox{\underline{$\gamma$}}:=(\gamma^1,\ldots,\gamma^{n_{\mbox{\tiny slip}}})\,.$$
Therefore, we get from (\ref{orientens2}) that
\begin{equation}\label{3rd-order-m}
p=\overline{m}\,\vectgam\,,\end{equation} where $\overline{m}$ is the third order tensor\footnote{The terminology ``tensor'' used here is just intended as ``matrix'' since $\overline{m}$ does not fulfill the rules of change of orthogonal bases for the last index and therefore is not a tensor in the usual sense.} defined as \begin{equation}\label{orientens3} \overline{m}_{ij\alpha}:=m^\alpha_{ij}=l^\alpha_i\nu^\alpha_j\quad\mbox{ for}\quad i,j=1,2,3\, \mbox{ and } \alpha=1,\ldots,n_{\mbox{\tiny slip}}\,.\end{equation}
Let $\mbox{\underline{$\eta$}}:=(\eta^1,\ldots,\eta^{n_{\mbox{\tiny slip}}})$ with $\eta^\alpha$ being a hardening variable in the $\alpha$-th slip system.
\subsection{The case with isotropic hardening}\label{micro-iso-crystal}

The starting point is the total energy 
\begin{equation}\label{total-energ1}
\mathcal{E}(u,\vectgam,\PL_p,\vecteta):=\int_\Omega\bigl[\Psi(\nabla u,\vectgam,\PL_p,\Curl\PL_p,\vecteta)-\la f,u\ra\bigr]\,dx\end{equation} where the free-energy density $\Psi$ is given in the additively separated form
\begin{eqnarray}\label{free-eng-crystal}
\Psi(\nabla u,\vectgam,\PL_p,\Curl\PL_p,\vecteta):
&=&\underbrace{\Psi^{\mbox{\scriptsize lin}}_e(\bvarepsilon_e)}_{\mbox{\small elastic energy}}\qquad
+\qquad\underbrace{\Psi^{\mbox{\scriptsize lin}}_{\mbox{\scriptsize
 micro}}(p,\PL_p)}_{\mbox{\small micro
energy}}\\
&& \nonumber +\,\,\underbrace{\Psi^{\mbox{\scriptsize lin}}_{\mbox{\scriptsize
curl}}(\Curl \PL_p)}_{\mbox{\small defect-like
energy (GND)}}\qquad+\hskip.4truecm\,\,\underbrace{\,\,\Psi_{\mbox{\scriptsize
 iso}}(\vecteta)}_{\begin{array}{c}
\mbox{\small hardening energy (SSD)}\end{array}}\,,
 \end{eqnarray} where
 \begin{eqnarray}\label{free-eng-expr-crystal} \nonumber\Psi^{\mbox{\scriptsize
lin}}_e(\bvarepsilon_e) &:=&\frac12\,\la\bvarepsilon_e,\C_{\mbox{\scriptsize iso}} \bvarepsilon_e\ra=\frac12\,\la\sym(\nabla u-p),\C_{\mbox{\scriptsize iso}}\sym(\nabla u-p)\ra\,, \\
\Psi^{\mbox{\scriptsize lin}}_{\mbox{\scriptsize
 micro}}(p,\PL_p)&:=&\frac12\,\mu\,H_\chi\,\norm{p-\PL_p}^2\,,\quad
 \Psi^{\mbox{\scriptsize lin}}_{\mbox{\scriptsize
curl}}(\Curl \PL_p)\,:=\,\frac12\,\mu\, L_c^2\norm{\Curl \PL_p}^2\,,\\ \nonumber\Psi_{\mbox{\scriptsize
 iso}}(\vecteta)&:=&\frac12\,\mu\,k_2\,\norm{\vecteta}^2\,=\,\frac12\mu\,k_2\,\sum_\alpha|\eta^\alpha|^2
 \,.\end{eqnarray}

Here,  $L_c\geq0$ is an energetic length scale which characterizes the contribution of the defect-like energy density to the system, $H_\chi$ is a positive nondimensional penalty constant, $k_2$  is a positive nondimensional isotropic hardening constant. \vskip.2truecm\noindent
The starting point for the derivation of the equations and inequalities describing the plasticity model is the two-field minimization formulation
\begin{equation}\label{micromorphic}
\mathcal{E}(u,\vectgam,\PL_p,\vecteta)\quad\rightarrow\quad \mbox{min. \,w.r.t. }\, (u,\PL_p)\,.\end{equation}
The first variations of  the total energy w.r.t. to the variables $u$ and $\PL_p$ lead to the balance equations in the next section.
\subsubsection{The balance equations}\label{Equi} The conventional macroscopic
force balance leads to the equation of equilibrium
\begin{equation}\label{EQL-crystal}
\div \sigma + f = 0
\end{equation}
in which $\sigma$ is the infinitesimal symmetric Cauchy stress and
$f$ is the body force.\\
An additional microscopic balance equation is obtained as follows. Precisely,  the first variation of the total energy w.r.t $\PL_p$ gives for every $\delta\PL_p\in C^\infty(\Omega,\mathbb{R}^{3\times 3})$, 
\begin{eqnarray}\label{microbal-weak1}\nonumber 0&=&\frac d{dt}\mathcal{E}(u,p,\PL_p+t\delta\PL_p,\vecteta)|_{t=0} \\
\nonumber &=&\int_\Omega\bigl[\mu\,H_\chi\,\la\PL_p-p,\delta\PL_p\ra 
 +\mu\,L^2_c,\la\Curl\PL_p,\Curl\delta\PL_p\ra\bigr]dx\\
\nonumber &=&     \int_\Omega\Bigl[\mu\,H_\chi\,\la\PL_p-p,\delta\PL_p\ra 
 +\mu\,L^2_c,\la\Curl\Curl\PL_p,\delta\PL_p\ra \\
\nonumber &&\hskip3truecm +     \sum_{i=1}^3\mbox{div}\Bigl(\mu
L_c^2\,\delta\PL_p^i\times(\Curl
\PL_p)^i\Bigr)\Bigr]dx\\             
 &=& \int_\Omega\la \mu\,H_\chi\,(\PL_p-p)
 +\mu\,L^2_c\,\la\Curl\Curl\PL_p,\delta\PL_p\ra \\
\nonumber &&\hskip3truecm +    \sum_{i=1}^3\int_{\partial\Omega}\mu
L_c^2\,\la\delta\PL_p^i\times(\Curl
\PL_p)^i,n\ra\,da,
\end{eqnarray} which on the one hand gives from the choice $\delta\PL_p\in C_c^\infty(\Omega,\mathbb{R}^{3\times 3})$  the micro-balance in strong formulation\footnote{Here, we have assumed uniform material constants for simplicity.}
\begin{equation}\label{EQL-micro-crystal}
\mu\,L^2_c\Curl\Curl\PL_p=-\,\mu\,H_\chi\,(\PL_p-p)\,.\end{equation} One the other hand we get
 \begin{equation}\label{EQL-micro-bc}  \sum_{i=1}^3\int_{\partial\Omega}\mu
L_c^2\,\la\delta\PL_p^i\times(\Curl
\PL_p)^i,n\ra\,da=0\quad\qquad\forall \delta\PL_p\in C^\infty(\Omega,\mathbb{R}^{3\times 3})\end{equation} which is satisfied
if we choose certain homogeneous  boundary conditions on the micro-distortion $\PL_p$. Following 
Gurtin \cite{GURT2004} and also Gurtin and Needleman \cite{GURTNEED2005}  we choose the simple boundary condition  \begin{equation}\label{bc-micro}
       \PL_p\times\,n|_{\Gamma_{\mbox{\scriptsize D}}}=0 \quad\mbox{ and }\quad (\Curl\PL_p)\times n|_{\partial\Omega\setminus\Gamma_{\mbox{\scriptsize D}}}=0
\end{equation}
which in the case of models in strain gradient plasticity, where $\PL_p$ is replaced by $p$ or $\bvarepsilon_p$ simply implies that there is no flow of the plastic distortion or plastic strain across the piece $\Gamma$ of the boundary $\partial\Omega$.

\subsubsection{The derivation of the dissipation inequality.}\label{deriv-diss-ineq-crystal}

The local free-energy imbalance states that
\begin{equation}
\dot{\Psi} - \la\sigma,\dot{\bvarepsilon}_e\ra - \la\sigma,\dot{p}\ra  \leq 0\
. \label{2ndlaw-crystal}
\end{equation}
Now we expand the first term, substitute (\ref{free-eng-crystal})-(\ref{free-eng-expr-crystal}) and get
\begin{equation}\label{exp-2ndlaw-crystal1}
\la\C_{\mbox{\scriptsize iso}}\,\bvarepsilon_e-\sigma,\dot{\bvarepsilon}_e\ra -\sum_\alpha\la\sigma,m^\alpha\ra\,\dot{\gamma}^\alpha+\sum_\alpha\frac{\partial \Psi^{\mbox{\scriptsize lin}}_{\mbox{\tiny
 micro}}}{\partial\gamma^\alpha}\,\dot{\gamma}^\alpha +\sum_\alpha\frac{\partial \Psi_{\mbox{\tiny
 iso}}}{\partial\eta^\alpha}\,\dot{\eta}^\alpha\leq0\,.\end{equation}
That is
\begin{equation}\label{exp-2ndlaw-crystal2}
\la\C_{\mbox{\scriptsize iso}}\,\bvarepsilon_e-\sigma,\dot{\bvarepsilon}_e\ra -\sum_\alpha\la\sigma,m^\alpha\ra\,\dot{\gamma}^\alpha-\sum_\alpha\mu\,H_\chi\,\la\PL_p-p,m^\alpha\ra\,\dot{\gamma}^\alpha +\mu\,k_2\,\sum_\alpha\eta^\alpha\,\dot{\eta}^\alpha\leq0\,.\end{equation} Therefore we obtain
\begin{equation}\label{exp-2ndlaw-crystal3}
 0\leq -\,\la\C_{\mbox{\scriptsize iso}}\,\bvarepsilon_e-\sigma,\dot{\bvarepsilon}_e\ra +\sum_\alpha\bigl[(\tau^\alpha +s^\alpha)\,\dot{\gamma}^\alpha +g^\alpha\,\dot{\eta}^\alpha\bigr]\\
\,=\, -\,\la\C_{\mbox{\scriptsize iso}}\,\bvarepsilon_e-\sigma,\dot{\bvarepsilon}_e\ra +\sum_\alpha\bigl[\tau^\alpha_{\mbox{\tiny E}}\,\dot{\gamma}^\alpha +g^\alpha\,\dot{\eta}^\alpha\bigr]\,,
\end{equation} 
where we set
\begin{eqnarray}
\tau^\alpha&:=&\la\sigma,m^\alpha\ra\quad\mbox{(resolved shear stress for the $\alpha$-th slip system)}\,,\\
s^\alpha&:=&\mu\,H_\chi\,\la\PL_p-p,m^\alpha\ra=-\mu\,L_c^2\,\la\Curl\Curl\PL_p,m^\alpha\ra\,,\\
g^\alpha&:=&-\,\mu\,k_2\,\eta^\alpha\quad\mbox{(thermodynamic force power-conjugate to $\dot{\eta}^\alpha$)} \,,\\
\tau^\alpha_{\mbox{\tiny E}}&:=&\tau^\alpha +s^\alpha=\la\Sigma_E,m^\alpha\ra\,,
\end{eqnarray} with  $\Sigma_E$  being the non-symmetric Eshelby-type stress tensor defined by
\begin{equation}\label{Eshelby-crystal} \Sigma_E:=\sigma+\mu\,H_\chi\,(\PL_p-p)\,\,=\,\,\sigma-\mu\,L_c^2\,\Curl\Curl\PL_p\,.
\end{equation}
Since the inequality (\ref{exp-2ndlaw-crystal3}) must be satisfied for whatever elastic-plastic deformation mechanism, inlcuding purely elastic ones (for which  $\dot{\gamma}^\alpha=0$, $\dot{\eta}^\alpha=0)$, equation  (\ref{exp-2ndlaw-crystal3}) implies the usual infinitesimal elastic stress-strain relation
\begin{eqnarray}
\nonumber \sigma = \C_{\mbox{\scriptsize iso}}\,\varepsilon_e&=&2\mu\, \sym(\nabla
u-p)+\lambda\, \tr(\nabla u-p)\id \\
&=& 2\mu\, (\sym(\nabla
u)-\bvarepsilon_p)+\lambda\, \tr(\nabla u)\id\label{elasticlaw}
\end{eqnarray}
and the local reduced dissipation inequality
\begin{equation}
\sum_\alpha\bigl[\tau_E^\alpha\,\dot{\gamma}^\alpha+g^\alpha\,\dot{\eta}^\alpha\bigr]\geq0
\end{equation}
which can also be written in compact form as 
\begin{equation}\label{reduced-diss-comp}
\sum_\alpha\la\Sigma_p^\alpha,\dot{\eta}_p^\alpha\ra\geq0\,,
\end{equation}where we define 
\begin{equation}\label{gen-var}\Sigma_p^\alpha:=(\tau^\alpha_E,g^\alpha)\quad\mbox{ and }\quad \Gamma_p^\alpha:=(\gamma^\alpha,\eta^\alpha)\,.
\end{equation}
\subsubsection{The flow rule}\label{flow-law}
We consider a yield function on the $\alpha$-th slip system
defined by
\begin{equation}
\label{yield-funct-kinematic} \phi(\Sigma^\alpha_p):= |\tau^\alpha_{\mbox{\tiny E}}|+g^\alpha - \yieldzero\quad\mbox{ for }\quad\Sigma^\alpha_p=(\tau^\alpha_{\mbox{\tiny E}},g^\alpha)\,.
\end{equation}
Here,  $\yieldzero$ is the yield stress of the material, that we assume to be constant on all slip systems and therefore, $\yieldlimit^\alpha:=\sigma_0-g^\alpha$ represents the current yield stress for the $\alpha$-th slip system{\footnote{Note that, for the sake of simplicity, the presented isotropic hardening
rule $g^\alpha$ does not involve latent hardening and the associated interaction matrix, see \cite{franciosizaoui91}
for a discussion on uniqueness in the presence of latent hardening.}. So the set
of admissible  generalized stresses for the $\alpha$-th slip system is defined as
\begin{equation}\label{admiss-stress-kin}
\mathcal{K}^\alpha:=\left\{\Sigma^\alpha_p=(\tau^\alpha_{\mbox{\tiny
E}},g^\alpha)\,\,|\,\,\phi(\Sigma_p^\alpha)
   \leq0,\,\,g^\alpha\leq0\right\}\,,\end{equation} with its interior $\mbox{Int}(\mathcal{K}^\alpha)$ and its boundary $\partial\mathcal{K}^\alpha$ being the generalized elastic region and the yield surface for the $\alpha$-th slip system, respectively.
   \vskip.1truecm\noindent
The  principle of maximum dissipation\footnote{The principle of maximum dissipation (PMD) is shown to be closely related to the so-called minimum principle  for the dissipation potential (MPDP) \cite{hacklfischer,hackl1997generalized,OrtRep99NEMD}, which states that the rate of the internal variables is the minimizer of a functional consisting of the sum of the rate of the free energy and the dissipation function with respect to appropriate boundary conditions. Notice that, as pointed out in \cite{EHN2016}, both PMD and MPDP are not physical principles but thermodynamically consistent selection rules which turn out to be convenient if no other information is available or if existing flow rules are to be extended to a more general situation.} associated with the $\alpha$-th slip system gives  us the normality
law\begin{equation}\label{normalcone} \dot{\Gamma}_p^\alpha\in
N_{\mathcal{K}^\alpha}(\Sigma^\alpha_p)\,,\end{equation}
where $\dsize N_{\mathcal{K}^\alpha}(\Sigma^\alpha_p)$
denotes the normal cone to $\mathcal{K}^\alpha$ at
$\dsize\Sigma_p^\alpha$. That is,  $\dot{\Gamma}^\alpha_p$ satisfies
\begin{equation}
\la\overline{\Sigma}^\alpha - \Sigma^\alpha_p,\dot{\Gamma}^\alpha_p\ra \leq
0\ \quad \mbox{for all}\ \overline{\Sigma}^\alpha \in {\mathcal{K}^\alpha}\ .
\label{normality2-alpha}
\end{equation} Notice that $N_{\mathcal{K}^\alpha}=\partial\Chi_{\mathcal{K}^\alpha}$, where
$\Chi_{\mathcal{K}^\alpha}$ denotes the indicator function of the set $\mathcal{K}^\alpha$
and $\partial\Chi_{\mathcal{K}^\alpha}$ denotes the subdifferential of the function
$\Chi_{\mathcal{K}^\alpha}$.\\ Whenever the yield surface $\partial\mathcal{K}^\alpha$ is
smooth at $\dsize\Sigma_p^\alpha$ then
$$\dot{\Gamma}^\alpha_p\in
N_{\mathcal{K}^\alpha}(\Sigma^\alpha_p)\quad\Rightarrow\quad\exists\lambda^\alpha\mbox{ such that }
\dot{\gamma}^\alpha=\lambda^\alpha\,\frac{\tau^\alpha_{\mbox{\tiny E}}}{|\tau^\alpha_{\mbox{\tiny E}}|}\quad\mbox{and}\quad\dot{\eta}^\alpha=\lambda^\alpha=|\dot{\gamma}^\alpha|$$ with the Karush-Kuhn
Tucker conditions: $\lambda^\alpha\geq0$, $\phi(\Sigma^\alpha_p)\leq0$ and $\lambda^\alpha\,\phi(\Sigma^\alpha_p)=0$\,.\\
Using convex analysis (Legendre-transformation) we find that
\begin{eqnarray}\label{dualflowlaw-crystal} &\underbrace{\dot{\Gamma}^\alpha_p\in
\partial\Chi_{\mathcal{K}^\alpha}(\Sigma^\alpha_p)}_{\mbox{\bf flow rule in its dual formulation for the $\alpha$-th slip system}}&\\\nonumber&&\\
\nonumber&\Updownarrow& \\\nonumber&&\\
&\underbrace{\Sigma^\alpha_p\in
\partial \Chi^*_{\mathcal{K}^\alpha} (\dot{\Gamma}^\alpha_p)\,}_{\mbox{\bf flow rule in its primal formulation  for the $\alpha$-th slip system}}&\label{primalflowlaw-crystal}
\end{eqnarray} where $\Chi^*_{\mathcal{K}^\alpha}$ is the Fenchel-Legendre dual of the function $\Chi_{\mathcal{K}^\alpha}$ denoted in this context by $\mathcal{D}^\alpha_{\mbox{\scriptsize iso}}$,
 the one-homogeneous dissipation function for the $\alpha$-th slip system. That is, for every $\Gamma^\alpha=(q^\alpha,\beta^\alpha)$,
\begin{eqnarray}\label{dissp-function-iso-slip}
\nonumber\mathcal{D}^\alpha_{\mbox{\scriptsize iso}}(\Gamma^\alpha)
&=&\sup\graffe{\la\Sigma^\alpha_p,\Gamma^\alpha\ra\,\,|\,\,\Sigma^\alpha_p\in\mathcal{K}^\alpha}\\
\nonumber&=&
\sup \graffe{\tau^\alpha_E\,q^\alpha +g^\alpha\beta^\alpha \,\,|\,\
\phi(\Sigma^\alpha_E,g^\alpha)\leq0,\,\,g^\alpha\leq0}\\&=&\left\{\begin{array}{ll}
      \yieldzero\,|q^\alpha| &\mbox{ if } |q^\alpha|\leq\beta^\alpha\,,\\
      \infty &\mbox{ otherwise.}\end{array}\right.
\end{eqnarray} We get from the definition of the subdifferential ($\Sigma^\alpha_p \in
\partial \Chi^*_{\mathcal{K}^\alpha} (\dot{\Gamma}^\alpha_p)$) that
\begin{equation}
\mathcal{D}^\alpha_{\mbox{\scriptsize iso}} (\Gamma^\alpha) \geq \mathcal{D}^\alpha_{\mbox{\scriptsize iso}}(\dot{\Gamma}^\alpha_p) +
\la\Sigma^\alpha_p,\Gamma^\alpha-\dot{\Gamma}^\alpha_p\ra\quad \mbox{for any
} \Gamma^\alpha.\label{dissineq-slip}\end{equation}
That is,
\begin{equation} \mathcal{D}^\alpha_{\mbox{\scriptsize iso}}(q^\alpha,\beta^\alpha)\geq \mathcal{D}^\alpha_{\mbox{\scriptsize iso}}(\dot{\gamma}^\alpha,\dot{\eta}^\alpha)+\tau^\alpha_{\mbox{\tiny E}}\,(q^\alpha-\dot{\gamma}^\alpha)+g^\alpha(\beta^\alpha-\dot{\eta}^\alpha)\quad \mbox{for any
} (q^\alpha,\beta^\alpha).\label{dissinequality-slip2}\end{equation}
\vskip.3truecm\noindent In the next sections, we present a complete mathematical analysis of the model including  both strong and weak formulations as well as a corresponding
existence result.
\subsubsection{Strong formulation
of the model}\label{strong-iso-crystal} To summarize, we have obtained the following strong
formulation for the microcurl model in the single crystal infinitesimal gradient plasticity case with
 isotropic  hardening.  Given $f\in \SFH^1(0,T;L^2(\Omega,\mathbb{R}^3))$, the goal is to find:
\begin{itemize}\item[(i)] the displacement $u\in \SFH^1(0,T; H^1_0(\Omega,{\Gamma_{\mbox{\scriptsize D}}},\mathbb{R}^3))$,
\item[(ii)] the infinitesimal plastic slips $\gamma^\alpha\in
\SFH^1(0,T;L^2(\Omega))$ for $\alpha=1,\ldots,n_{\mbox{\tiny slip}}$, the infinitesimal micro-distortion $\PL_p\in\SFH^1(0,T;H(\Curl;\Omega, \mathbb{R}^{3\times 3})$, with $ \Curl\Curl\PL_p \in \SFH^1(0,T;
L^2(\Omega,\BBR^{3\times 3}))$\end{itemize}
 such that the content of Table \ref{table:micro-isohard-crystal} holds.
 
 \begin{table}[ht!]{\footnotesize\begin{center}
\begin{tabular}{|ll|}\hline &\qquad\\
 {\em Additive split of distortion:}& $\nabla u =e +p$,\quad $\bvarepsilon_e=\mbox{sym}\,e$,\quad $\bvarepsilon_p=\sym p$\\
{\em Plastic distortion in slip system:}&$p=\dsize\sum_{\alpha=1}^{n_{\mbox{\tiny silp}}}      \gamma^\alpha\,m^\alpha$ with $m^\alpha= l^\alpha\otimes \nu^\alpha$, \quad $\tr(p)=0$\\
{\em Equilibrium:} & $\mbox{Div}\,\sigma +f=0$ with
$\sigma=\C_{\mbox{\tiny iso}}\bvarepsilon_e=\C_{\mbox{\tiny iso}}(\sym\nabla u-\bvarepsilon_p)$\\
{\em Microbalance:} & $\mu\,L_c^2\,\Curl\Curl\PL_p=-\mu\,H_\chi\,(\PL_p-p)$,\\&\\
 {\em Free energy:} &
$\frac12\,\langle\C_{\mbox{\tiny iso}}\bvarepsilon^e,\bvarepsilon^e\rangle+\frac12\,\mu\,H_\chi\,\norm{p-\PL_p}^2$\\
&\qquad\quad$+\,\frac12\,\mu\,
L^2_c\,\norm{\Curl \PL_p}^2+\frac12\,\mu\, k_2\,\sum_\alpha|\eta^\alpha|^2$\\&\\
{\em Yield condition in $\alpha$-th slip system:} &
$|\tau^\alpha_{\mbox{\tiny E}}|+g^\alpha-\yieldzero\leq0$\\
 {\em where } & $\tau^\alpha_{\mbox{\tiny E}}:=\la\Sigma_E,m^\alpha\ra$ with\\
 & $\Sigma_E:=\sigma+\mu\,H_\chi\,(\PL_p-p)=\sigma-\mu\,L_c^2\,\Curl\Curl\PL_p$\\
 & $g^\alpha=-\mu\,k_2\,\eta^\alpha$
 \\&\\{\em  Dissipation inequality in $\alpha$-th slip system:} &
 $\tau^\alpha_{\mbox{\tiny E}}\,\dot{\gamma}^\alpha +g^\alpha\,\dot{\eta}^\alpha\geq0$\\
 {\em Dissipation function in $\alpha$-th slip system:} &$\mathcal{D}^\alpha_{\mbox{\tiny iso}}(q^\alpha,\beta^\alpha):=\left\{\begin{array}{ll}\yieldzero |q^\alpha| &\mbox{ if } |q^\alpha|\leq\beta^\alpha,\\ \infty &\mbox{ otherwise}\end{array}\right.$\\&\\
 {\em Flow rules in primal form:} &
 $(\tau^\alpha_{\mbox{\tiny E}},g^\alpha)\in\partial \mathcal{D}^\alpha_{\mbox{\scriptsize iso}}(\dot{\gamma}^\alpha,\dot{\eta}^\alpha)$\\&\\
{\em Flow rules in dual form:}
&$\dot{\gamma}^\alpha=\lambda^\alpha\,\dsize\frac{\tau^\alpha_{\mbox{\tiny E}}}{|\tau^\alpha_{\mbox{\tiny E}}|},\quad\qquad \dot{\eta}^\alpha=\lambda^\alpha=|\dot{\gamma}^\alpha|$\\&\\
{\em KKT conditions:} &$\lambda^\alpha\geq0$, \quad $\phi(\tau^\alpha_E,g^\alpha)\leq0$,
\quad $\lambda^\alpha\,\phi(\tau^\alpha_{\mbox{\tiny E}},g^\alpha)=0$\\
 {\em Boundary conditions for $\PL_p$:} & $\PL_p\times{n}=0$ on
 $\Gamma_{\mbox{\scriptsize D}}$,\,\, $(\Curl \PL_p)\times{n}=0$ on $\partial\Omega\setminus\Gamma_{\mbox{\scriptsize D}}$\\
 {\em Function space for $\PL_p$:} & $\PL_p(t,\cdot)\in \mbox{H}(\mbox{Curl};\;\Omega,\,\BBR^{3\times 3})$\\
 \hline
\end{tabular}\caption{\footnotesize The microcurl model in single crystal gradient plasticity  with isotropic hardening. }\label{table:micro-isohard-crystal}\end{center}}\end{table}
 
\subsubsection{Weak formulation of the model}\label{wf-iso-crystal} Assume that the problem in Section \ref{strong-iso-crystal} has a solution
$(u,\vectgam,\PL_p,\vecteta)$.  We will extensively make use of the identity (\ref{3rd-order-m}). Let $v\in H^1(\Omega,\mathbb{R}^3)$
with $v_{|\Gamma_D}=0$. Multiply the equilibrium equation with
$v-\dot{u}$ and integrate in space by parts and use the
symmetry of $\sigma$ and the elasticity relation to get
\begin{equation}\label{weak-eq1}
\int_{\Omega}\la\C_{\mbox{\scriptsize iso}}\sym(\nabla
u-\overline{m}\,\underline{\gamma}),\mbox{sym}(\nabla v-\nabla\dot{u})\ra\, dx=\int_\Omega
f(v-\dot{u})\,dx\, .
\end{equation}
Now,
for any $X\in C^\infty(\overline{\Omega},\sL(3))$ such that
$X\times\,n=0$ on $\Gamma$ 
we integrate (\ref{EQL-micro-crystal}) over $\Omega$, integrate by
parts the term with Curl\,Curl using the boundary conditions
$$(X-\dot{\PL}_p)\times\,n=0\mbox{ on }\Gamma_{\mbox{\scriptsize D}},\qquad
\mbox{Curl}\,\PL_p\times\,n=0\mbox{ on }
\partial\Omega\setminus\Gamma_{\mbox{\scriptsize D}}$$ and get
\begin{equation}\label{EQL-microweak}
\int_\Omega\Bigl[\mu\,L_c^2\la\Curl\PL_p,\Curl X-\Curl\dot{\PL}_p\ra+\mu\,H_\chi\,\la\PL_p-\overline{m}\,\underline{\gamma},X-\dot{\PL}_p\ra\Bigr]\,dx=0\,.
\end{equation}
Moreover,
for any $\underline{q}=(q^1,\ldots,q^{n_{\mbox{\tiny slip}}})$ with $q^\alpha\in C^\infty(\overline{\Omega})$  and any $\underline{\beta}=(\beta^1,\ldots,\beta^{n_{\mbox{\tiny slip}}})$ with $\beta^\alpha\in L^2(\Omega)$, summing (\ref{dissinequality-slip2}) over $\alpha=1,\ldots,n_{\mbox{\tiny slip}}$ and  integrating over $\Omega$, we  get

\begin{eqnarray}\label{dissinequality2-wk}
\nonumber &&\hskip-1truecm \int_\Omega\mathcal{D}_{\mbox{\tiny iso}}(\underline{q},\underline{\beta})\,dx
-\int_\Omega\mathcal{D}_{\mbox{\tiny iso}}(\dot{\vectgam},\dot{\vecteta})\,dx -\int_\Omega
\la\C_{\mbox{\tiny iso}}\sym(\nabla u-\overline{m}\,\underline{\gamma}),\sym(\overline{m}\,\underline{q}-\overline{m}\,\dot{\underline{\gamma}})\ra\,dx\\
&&\quad+\int_\Omega\Bigl[-\mu\,H_\chi\,\la\PL_p-\overline{m}\,\underline{\gamma},\overline{m}\,\underline{q}-\overline{m}\,\dot{\underline{\gamma}}\ra
+\mu\,k_2\,\la\vecteta,\underline{\beta}-\dot{\vecteta}\ra\,\Bigr]dx\geq0\,.\end{eqnarray}
where
\begin{equation}\label{overall-diss}
\mathcal{D}_{\mbox{\tiny iso}}(\underline{q},\underline{\beta}):=\sum_\alpha \mathcal{D}^\alpha_{\mbox{\tiny iso}}(q^\alpha,\beta^\alpha)\,.\end{equation}

Now adding up (\ref{weak-eq1})-(\ref{dissinequality2-wk}) we
get the following weak formulation of the problem set in Section
\ref{strong-iso-crystal} in the form of a variational inequality:
\begin{eqnarray}\label{weak-form-iso-crystal}
\nonumber&&\hskip-.7truecm\int_\Omega\Bigl[\la\C_{\mbox{\tiny iso}}\sym(\nabla
u-\overline{m}\,\underline{\gamma}),\sym(\nabla v-\overline{m}\,\underline{q})-\sym(\nabla\dot{u}-\overline{m}\,\dot{\underline{\gamma}})\ra+ \mu\,L_c^2\la\Curl\PL_p,\Curl X-\Curl\dot{\PL}_p\ra\\
\nonumber&&\hskip1truecm+\,\,\mu\,H_\chi\,\la\PL_p-\overline{m}\,\underline{\gamma}, (X-\overline{m}\,\underline{q})-(\dot{\PL}_p-\overline{m}\,\dot{\underline{\gamma}})\ra+\mu\,k_2\,\la\vecteta,\underline{\beta}-\dot{\vecteta}\ra\,\Bigr]dx\\
&&\hskip3truecm  +
\int_\Omega\mathcal{D}_{\mbox{\tiny iso}}(\underline{q},\underline{\beta})\,dx
-\int_\Omega\mathcal{D}_{\mbox{\tiny iso}}(\dot{\vectgam},\dot{\vecteta})\,dx\,\geq\, \int_\Omega
f\,(v-\dot{u})\,dx\,.\end{eqnarray}

\subsubsection{Existence result for the weak formulation}\label{exist-crystal-iso}
To prove the existence result for the weak formulation
(\ref{weak-form-iso-crystal}), we closely follow the abstract machinery developed by
Han and Reddy in \cite{Han-ReddyBook} for mathematical problems in geometrically linear
classical plasticity and used for instance in \cite{DEMR1, REM, NCA,
EBONEFF, ENR2015} for models of gradient plasticity. To this aim, equation 
(\ref{weak-form-iso-crystal}) is written as the variational inequality of
the second kind: find $w=(u,\vectgam,\PL_p,\vecteta)\in \SFH^1(0,T;\SFZ)$
such that $w(0)=0$,  $\dot{w}(t)\in \SFW$ for a.e. $t\in[0,T]$ and
\begin{equation}\label{wf}
\ba(w,z-\dot{w})+j(z)-j(\dot{w})\geq \langle
\ell,z-\dot{w}\rangle\mbox{ for every } z\in \SFW\mbox{ and for a.e.
}t\in[0,T]\,,\end{equation} where $\SFZ$ is a suitable Hilbert space
and $\SFW$ is some closed, convex subset of $\SFZ$ to be constructed
later,
\begin{eqnarray}
\ba(w,z)&=&\int_\Omega\Bigl[\la\C_{\mbox{\scriptsize iso}}\sym(\nabla
u-\overline{m}\,\underline{\gamma}),\sym(\nabla v-\overline{m}\,\underline{q})\ra+ \mu\,L_c^2\la\Curl\PL_p,\Curl X\ra\\
\nonumber&&\hskip2truecm+\,\,\mu\,H_\chi\,\la\PL_p-\overline{m}\,\underline{\gamma}, X-\overline{m}\,\underline{q}\ra+\mu\,k_2\,\la\vecteta,\underline{\beta}\ra\,\Bigr]dx\,,
\label{bilin-iso-spin}
\\\nonumber\\
j(z)&=&\int_\Omega\mathcal{D}_{\mbox{\tiny iso}}(\underline{q},\underline{\beta})\,dx\,,\label{functional-isospin}\\
\langle \ell,z\rangle&=&\int_\Omega
f\,v\,dx\,,\label{lin-form}\end{eqnarray} for
$w=(u,\vectgam,\PL_p,\vecteta)$ and $z=(v,\underline{q},X,\underline{\beta})$ in
$\SFZ$.
\vskip.2truecm\noindent
The Hilbert space $\SFZ$ and the closed convex subset $\SFW$ are
constructed in such a way that the functionals $\ba$, $j$ and
$\ell$ satisfy the assumptions in the abstract result in
\cite[Theorem 6.19]{Han-ReddyBook}. The key issue here is the
coercivity of the bilinear form $\ba$ on the set $\SFW$, that is, $a(z,z)\geq C\norm{z}^2_Z$ for every $z\in\SFW$ and for some $C>0$.
\vskip.2truecm\noindent
We let
\begin{eqnarray}
\SFV&:=&\mathsf{H}^1_0(\Omega,{\Gamma_{\mbox{\scriptsize
D}}},\mathbb{R}^3)=\bigl\{v\in \mathsf{H}^1(\Omega,\mathbb{R}^3)\,\,|\,\,
v_{|\Gamma_{\mbox{\scriptsize
D}}}=0\bigr\}\,,\label{space-v-crytsal}\\
 \SFP&:=&\mbox{L}^2(\Omega,\,\mathbb{R}^{n_{\mbox{\tiny slip}}})\,,\label{space-epsp-crystal}\\
\SFQ&:=&\mbox{H}_0(\mbox{Curl};\,\Omega,\,\Gamma_{\mbox{\scriptsize D}},\sL(3))\,,\label{space-p-crystal}\\
\Lambda&:=& L^2(\Omega,\,\mathbb{R}^{n_{\mbox{\tiny slip}}})\,,\label{space-hardvar-crystal}\\
   \SFZ&:=&\SFV\times \SFP\times\SFQ\times\Lambda\,,\label{product-space-crystal}\\
   \SFW&:=&\Bigl\{z=(v,\underline{q},X,\underline{\beta})\in\SFZ\,\,\,|\,\,\,|q^\alpha|\leq\beta^\alpha\,,\quad \alpha=1,\ldots,n_{\mbox{\tiny slip}}\Bigr\}\,,\label{set-W-crystal}
\end{eqnarray} and define the norms
\begin{eqnarray}
\nonumber&&  \norm{v}_V:=\norm{\nabla v}_{L^2}\,,\label{norm-V} \quad \norm{\underline{q}}^2_P=\dsize\sum_\alpha\norm{q^\alpha}^2_{L^2}\,,\quad  \norm{\underline{\beta}}^2_{\Lambda}=\dsize\sum_\alpha\norm{\beta^\alpha}^2_{L^2}\,,\quad
\norm{X}_{Q}:=\norm{X}_{\mbox{\scriptsize H}(\mbox{\scriptsize
Curl};\Omega)},\label{norm-Q}\\
&&\norm{z}^2_{Z}:=\norm{v}^2_{V} +\norm{\underline{q}}^2_{P}+\norm{X}^2_{Q}
+\norm{\underline{\beta}}^2_{\Lambda}\quad\mbox{ for }
z=(v,\underline{q},X,\underline{\beta})\in \SFZ\,. \label{norm-Z}
\end{eqnarray}

Let us show that the bilinear form $\ba$ is coercive on $\SFW$.
 Let therefore $z=(v,\underline{q},X,\underline{\beta})\in \SFW$. First of all notice that 
 \begin{equation}\label{normP-dom}
 \norm{\overline{m}\,\underline{q}}_{L^2}\leq\norm{\underline{q}}_P\,\leq\,\norm{\underline{\beta}}_{\Lambda}\,.\end{equation}
 So,
\begin{eqnarray}
\nonumber\ba(z,z)& \geq&  m_0\norm{\sym(\nabla v-\overline{m}\,\underline{q})}^2_{L^2}\mbox{
(from (\ref{ellipticityC}))} +\, \mu\,H_\chi\,\norm{X-\overline{m}\,\underline{q}}^2_{L^2}+\,\mu\,
L_c^2\norm{\Curl X}_{L^2}^2 \\
\nonumber &&\hskip.5truecm+\mu\, k_2\norm{\underline{\beta}}_{\Lambda}^2 \\
\nonumber&=&m_0\bigl[\norm{\sym\nabla v}^2_{L^2} +\norm{\sym(\overline{m}\,\underline{q})}^2_2
-2\la\sym\nabla v, \sym(\overline{m}\,\underline{q})\ra_{L^2}\bigr]+\mu\,H_\chi\,\bigl[\norm{X}^2_{L^2} \\
\nonumber &&\hskip.5truecm+\,\norm{\overline{m}\,\underline{q}}^2_{L^2}-\,2\la X, \overline{m}\,\underline{q}\ra_{L^2}\bigr]  +\,\mu L_c^2\,\norm{\Curl X}_{L^2}^2
+\mu\,k_2\norm{\underline{\beta}}_{\Lambda}^2\\
\nonumber&\geq &m_0\left[\norm{\sym\nabla v}^2_{L^2} +\norm{\sym(\overline{m}\,\underline{q})}^2_{L^2} -\theta\norm{\sym(\nabla
v)}_{L^2}^2-\frac1\theta\norm{\sym(\overline{m}\,\underline{q})}_{L^2}^2\right]\\
\nonumber&& \hskip.5truecm +\, \mu\,H_\chi\,\left[\norm{X}^2_{L^2} +\norm{
\overline{m}\,\underline{q}}^2_{L^2} -\theta\norm{X}_{L^2}^2-\frac1\theta\norm{\overline{m}\,\underline{q}}_{L^2}^2\right]\mbox{ (Young's inequality)}\\
\nonumber&& \hskip.5truecm+\,\mu \,L_c^2\norm{\Curl X}_{L^2}^2+\frac12\,\mu\, k_2\norm{\underline{\beta}}_{\Lambda}^2  +\,\frac12\,\mu\, k_2\,\norm{\underline{q}}_{P}^2\quad\mbox{ \,\,(using second $\leq$ in (\ref{normP-dom}))} \\
\nonumber &=& m_0(1-\theta)\norm{\sym\nabla
v}^2_{L^2} +m_0\Bigl(1-\frac1\theta\Bigr)\norm{\sym(\overline{m}\,\underline{q})}_{L^2}^2+\mu\,H_\chi\Bigl(1-\frac1\theta\Bigr)\norm{\underline{q}}^2_P\\
\nonumber &&\hskip.5truecm +\,\frac12\,\mu\,
k_2\,\norm{\underline{q}}_{P}^2+\mu\,H_\chi\,(1-\theta)\norm{X}_{L^2}^2 +\mu \,L_c^2\norm{\Curl X}_2^2 +
\frac12\,\mu\,k_2\,\norm{\underline{\beta}}_{\Lambda}^2\\
 &\geq& m_0(1-\theta)\norm{\sym\nabla
v}^2_{L^2}+\left[(m_0+\mu\,H_\chi)\Bigl(1-\frac1\theta\Bigr)+\frac12\,\mu\,
k_2\right]\norm{\underline{q}}_{P}^2\\
\nonumber && \hskip.5truecm+\,\mu\,H_\chi\,(1-\theta)\norm{X}_{L^2}^2 +\mu \,L_c^2\norm{\Curl X}_2^2 +
\frac12\,\mu\,k_2\,\norm{\underline{\beta}}_{\Lambda}^2\,\quad\mbox{ (for $0<\theta<1$)\,.}\label{coercive-a}\end{eqnarray}

So, since the hardening constant $k_2>0$, it is possible to choose $\theta$ such that $$\displaystyle\frac{m_0+\mu\,H_\chi\,}{m_0+\mu\,H_\chi+\frac
12\,\mu\,k_2}< \theta<1,$$  we always are able to  find some constant
$C(\theta,m_0,\mu,H_\chi,k_2,L_c,\Omega)>0$ such that
\begin{equation}\label{coerc-iso-crystal}a(z,z)\geq C\left[\norm{v}_V^2+\norm{\underline{q}}^2_{P}+\norm{X}^2_{\mbox{\scriptsize H}(\mbox{\scriptsize
Curl};\Omega)}+\norm{\underline{\beta}}_{\Lambda}^2
\right]=C\norm{z}^2_{Z}\quad\forall z=(v,\underline{q},X,\underline{\beta})\in \SFW\,.\end{equation} This shows existence for the microcurl model in single gradient plasticity with isotropic hardening.
\subsubsection{Uniqueness of the weak/strong solution}\label{uniqueness} As shown in \cite{EHN2016} for a  canonical rate-independent model of
geometrically linear isotropic gradient plasticity with isotropic
hardening and plastic spin, the uniqueness of the solution for our model can  be obtained  similarly. To this aim, notice that if $(u,\vectgam, \PL_p,\vecteta)$ is a weak solution of the model, then $(u,\vectgam,\PL_p,\vecteta)$ is also a strong solution. In fact, choosing appropriatetly test functions in the variational inequality (\ref{weak-form-iso-crystal}), we obtain  both equilibrium and microbalance equations on the one hand. The latter which is 
$$\mu\,L^2_c\,\Curl\Curl\PL_p=-\mu\,H_\chi\,(\PL_p-p)$$ is satisfied first in the distributional sense and hence is satisfied also in the $L^2$-sense since the right hand side is in $L^2(\Omega,\mathbb{R}^{3\times 3})$. Therefore, it follows that  $\Curl\Curl\PL_p$ is also in $L^2(\Omega,\mathbb{R}^{3\times 3})$. Now going back to (\ref{weak-form-iso-crystal}), we also derive the boundary condition $(\Curl \PL_p)\times n|_{\partial\Omega\setminus\Gamma_{\mbox{\scriptsize D}}}=0$ which is now justified because we derived that $\Curl\PL_p\in \mbox{H}(\mbox{Curl};\,\Omega,\,\sL(3))$.\\
On the other hand, we also obtain from (\ref{weak-form-iso-crystal})  the following set of  inequatlities
\begin{equation}\label{flow-unique1}
\la\dot{\Gamma}^\alpha_p,\Sigma^\alpha-\Sigma^\alpha_p\ra\leq0\quad\forall\Sigma^\alpha\in\mathcal{K}^\alpha,\quad\alpha=1,\ldots,n_{\mbox{\tiny silp}}\,\end{equation} and hence, 
$(u,\vectgam,\PL_p,\vecteta)$ is a strong solution.

Now let us consider two solutions $w_i:=(u_i,\underline{\rgam_i}, \PL_{p\,i},\underline{\reta_i})$, $i=1,2$ of  (\ref{weak-form-iso-crystal}) satisfying the same initial conditions, let  $\Gamma_{p_{\,i}}^\alpha=(\gamma_i^\alpha,\eta^\alpha_i)$ and $\Sigma^\alpha_{p\,i}:=(\tau^\alpha_{E_{\,i}},g_i^\alpha)$ be the corresponding stresses. That is,
\begin{eqnarray}
\tau^\alpha_{E_{\,i}}&=&\la\Sigma_{E_{\,i}},m^\alpha\ra=\la\sigma_i +\mu\,H_\chi\,(\PL_{p_{\,i}}-p_i),m^\alpha\ra=\la\sigma_i-\mu\,H_\chi\Curl\Curl\PL_{p_{\,i}},m^\alpha\ra\,,\\
g^\alpha_i&=&-\mu\,k_2\,\eta^\alpha_i\,\end{eqnarray}
so that $\Gamma^\alpha_{p_{\,i}}$ and $\Sigma^\alpha_{p_{\,i}}$ satisfy

\begin{equation}\label{flow-unique2}
\la\dot{\Gamma}^\alpha_{p_{\,1}},\Sigma^\alpha-\Sigma^\alpha_{p_{\,1}}\ra\leq0\quad\mbox{ and }\quad
\la\dot{\Gamma}^\alpha_{p_{\,2}},\Sigma^\alpha-\Sigma^\alpha_{p_{\,2}}\ra\leq0 \quad\forall\Sigma^\alpha\in\mathcal{K}^\alpha\,.\end{equation}

Now choose $\Sigma^\alpha=\Sigma^\alpha_{p_{\,2}}$ in
(\ref{flow-unique2})$_1$ and $\Sigma^\alpha=\Sigma^\alpha_{p_{\,1}}$ in
(\ref{flow-unique2})$_2$ and add up to get
\begin{equation}\label{flow-unique3}
\la\Sigma^\alpha_{p_{\,2}}-\Sigma^\alpha_{p_{\,1}},\dot{\Gamma}^\alpha_{p_{\,1}}-\dot{\Gamma}^\alpha_{p_{\,2}}\ra\leq0\,.\end{equation}
That is,
\begin{eqnarray}\label{flow-unique4}
&&\nonumber\la\sigma_2-\sigma_1,m^\alpha\,\dot{\gamma}_1^\alpha-m^\alpha\,\dot{\gamma}_2^\alpha\ra+\mu\,H_\chi\,\la(\PL_{p_{\,2}}-\PL_{p_{\,1}})-(p_2-p_1),m^\alpha\,\dot{\gamma}^\alpha_1-m^\alpha\,\dot{\gamma}^\alpha_2\ra\\
&&\hskip4truecm+\,(g_2^\alpha-g_1^\alpha)(\dot{\eta}^\alpha_1-\dot{\eta}^\alpha_2)\leq0\end{eqnarray}
and adding up over $\alpha$, we get

\begin{equation}\label{flow-unique4}
\la\sigma_2-\sigma_1,\dot{p}_1-\dot{p}_2\ra+\mu\,H_\chi\,\la(\PL_{p_{\,2}}-\PL_{p_{\,1}})-(p_2-p_1),\dot{p}_1-\dot{p}_2\ra+\la\underline{g_2}-\underline{g_1},\underline{\dot{\eta}_1}-\underline{\dot{\eta}_2}\ra\leq0\end{equation}
 Now, substitute $\sym p_i=\sym\nabla u_i-\C_{\mbox{\tiny iso}}^{-1}\sigma_i$
obtained from the elasticity relation, into the expression $\la\sigma_2-\sigma_1,\dot{p}_1-\dot{p}_2\ra$ and $\dot{p}_i=\dot{\PL}_{_{p\,i}}+\dsize\frac{L_c^2}{H_\chi}\,\Curl\Curl\dot{\PL}_{p_{\,i}}$ obtained from the microbalance equation into the expression  $\la\PL_{p_{\,2}}-\PL_{p_{\,1}},\dot{p}_1-\dot{p}_2\ra$ and get from (\ref{flow-unique4}) that

\begin{eqnarray}\label{flow-unique5}
&&\nonumber\hskip-2.3truecm\la\sigma_2-\sigma_1,\C_{\mbox{\tiny iso}}^{-1}(\dot{\sigma}_2-\dot{\sigma}_1)\ra +\mu\,H_\chi\la p_1-p_2,\dot{p}_1-\dot{p}_2\ra +\mu\,H_\chi\la\PL_{p_{\,2}}-\PL_{p_{\,1}},\dot{\PL}_{p_{\,1}}-\dot{\PL}_{p_{\,2}}\ra\\
&&\hskip.1truecm +\,\mu\,L_c^2\la\PL_{p_{\,2}}-\PL_{p_{\,1}},\Curl\Curl\PL_{p_{\,1}}-\Curl\Curl\PL_{p_{\,2}}\ra +\la\underline{g_2}-\underline{g_1},\underline{\dot{\eta}_1}-\underline{\dot{\eta}_2}\ra\\
&&\nonumber\hskip4truecm\leq\, \la\sigma_1-\sigma_2,\sym(\nabla u_1)-\sym(\nabla u_2)\ra\end{eqnarray}
Now for every $t\in[0,T]$, we integrate (\ref{flow-unique5}) over $\Omega\times(0,t)$ using the boundary conditions on $\PL_{p_{\,i}}$ and using the fact that
$$\int_\Omega\la\sigma_1-\sigma_2,\sym\nabla u_1-\sym\nabla u_2\ra\,dx=0\,,$$ we get
\begin{eqnarray}\label{flow-unique6}
&&\nonumber\hskip-.8truecm\int_0^t\frac{d}{ds}\Bigl[\norm{\C_{\mbox{\tiny iso}}^{-1/2}(\sigma_2(s)-\sigma_1(s))}^2_{L^2}+\mu\,H_\chi\,\norm{p_1(s)-p_2(s)}^2_{L^2}-\,\,\mu\,H_\chi\,\norm{\PL_{p_{\,2}}(s)-\PL_{p_{\,1}}(s)}^2_{L^2}\\
&&\nonumber\hskip1truecm-\,\,\mu\,L_c^2\,\norm{\Curl\PL_{p_{\,2}}(s)-\Curl\PL_{p_{\,1}}(s)}^2_{L^2}+\mu\,k_2\,\norm{\underline{{\eta}_1}(s)-\underline{{\eta}_2}(s)}^2_{L^2}\Bigr]\,ds\leq0\,.
\end{eqnarray}
Therefore, we obtain
\begin{eqnarray}\label{flow-unique7}
&&\nonumber\hskip-2truecm\norm{\C_{\mbox{\tiny iso}}^{-1/2}(\sigma_2-\sigma_1)}^2_{L^2}+\mu\,H_\chi\,\norm{p_1-p_2}^2_{L^2}+\mu\,k_2\,\norm{\underline{{\eta}_1}-\underline{{\eta}_2}}^2_{L^2}\\
&&\hskip1truecm\leq\,\mu\,H_\chi\,\norm{\PL_{p_{\,2}}-\PL_{p_{\,1}}}^2_{L^2}+\mu\,L_c^2\,\norm{\Curl\PL_{p_{\,2}}-\Curl\PL_{p_{\,1}}}^2_{L^2}.
\end{eqnarray}
On the other hand, we write the micro-balance equation for $p_i$ and $\PL_{p_{\,i}}$ with $i=1,2$, as 
\begin{eqnarray}\label{flow-unique8-1}\mu\,H_\chi\,p_1&=&\mu\,L_c^2\,\Curl\Curl\PL_{p_{\,1}}+\mu\,H_\chi\PL_{p_{\,1}}\,,\\
\label{flow-unique8-2}\mu\,H_\chi\,p_2&=&\mu\,L_c^2\,\Curl\Curl\PL_{p_{\,2}}+\mu\,H_\chi\PL_{p_{\,2}}\,,\end{eqnarray}
then we subtract, take the scalar product with $\PL_{p_{\,1}}-\PL_{p_{\,2}}$, integrate using the boundary condition
$$(\PL_{p_{\,1}}-\PL_{p_{\,2}})\times n|_{\Gamma_{\mbox{\scriptsize D}}}=0\quad\mbox{ and }\quad(\Curl\PL_{p_{\,1}}-\Curl\PL_{p_{\,2}})\times n|_{\partial\Omega\setminus\Gamma_{\mbox{\scriptsize D}}}=0\,$$
 and get
\begin{equation}\label{flow-unique9}
\mu\,H_\chi\int_\Omega\la p_1-p_2,\PL_{p_{\,1}}-\PL_{p_{\,2}}\ra\,dx=\mu\,L_c^2\norm{\Curl\PL_{p_{\,1}}-\Curl\PL_{p_{\,2}}}^2_{L^2}+\mu\,H_\chi\norm{\PL_{p_{\,1}}-\PL_{p_{\,2}}}^2_{L^2}.\end{equation}
Therefore, we obtain from (\ref{flow-unique7}) and (\ref{flow-unique9}) that 
\begin{eqnarray}\label{flow-unique10}
\mu\,H_\chi\,\norm{p_1-p_2}^2_{L^2}&\leq& \mu\,L_c^2\norm{\Curl\PL_{p_{\,1}}-\Curl\PL_{p_{\,2}}}^2_{L^2}+\mu\,H_\chi\norm{\PL_{p_{\,1}}-\PL_{p_{\,2}}}^2_{L^2}\\
&\leq&\mu\,H_\chi\norm{p_1-p_2}_{L^2}\norm{\PL_{p_{\,1}}-\PL_{p_{\,2}}}_{L^2}
\end{eqnarray}
which implies that 
\begin{equation}\label{flow-unique11}
\norm{p_1-p_2}_{L^2}=\norm{\PL_{p_{\,1}}-\PL_{p_{\,2}}}_{L^2}\quad\mbox{ and hence, }\quad
\norm{\Curl\PL_{p_{\,1}}-\Curl\PL_{p_{\,2}}}_{L^2}=0\,.\end{equation}
Now, going back to (\ref{flow-unique7}), we get 
\begin{equation}\label{flow-unique12}
\norm{\C_{\mbox{\tiny iso}}^{-1/2}(\sigma_2-\sigma_1)}^2_{L^2}+\mu\,k_2\,\norm{\underline{{\eta}_1}-\underline{{\eta}_2}}^2_{L^2}\leq0\,.\end{equation}
Hence, we obtain so far,
$$\sigma_1-\sigma_2,\quad \underline{{\eta}_1}=\underline{{\eta}_2}\quad\mbox{ and }
\quad\Curl\PL_{p_{\,1}}=\Curl\PL_{p_{\,2}}\quad\Rightarrow\quad\tau^\alpha_{E_{\,1}}=\tau^\alpha_{E_{\,2}}\,.$$

Now, let us prove that $\gamma^\alpha_1=\gamma^\alpha_2$ for every $\alpha$. In fact, from the definition of
the normal cone it follows that $\dot{\Gamma}^\alpha_{p_{\,i}}=0$ that is, $\dot{\gamma}^\alpha_i=\dot{\eta}^\alpha_i=0$
 inside the elastic
domain Int$(\mathcal{K}^\alpha)$ (for the $\alpha$-slip system), which from the initial conditions imply that
$\gamma^\alpha_i=0$ inside Int$(\mathcal{K}^\alpha)$. Now, looking at the flow rule in dual
form (for the $\alpha$-slip system) in Table \ref{table:micro-isohard-crystal}, we obtain from $\tau^\alpha_{E_{\,1}}=\tau^\alpha_{E_{\,2}}$ that $\dot{\gamma}^\alpha_1=\dot{\gamma}^\alpha_2$ which implies
that $\gamma^\alpha_1=\gamma^\alpha_2$ from the initial conditions. Therefore, we obtain $p_1=p_2$ which implies from 
 (\ref{flow-unique11})$_1$ that $\PL_{p_{\,1}}=\PL_{p_{\,2}}$.
 
 Now, it remains to show that $u_1=u_2$.  This is obtained exactly as in \cite{EHN2016}. We repeat the proof here just for the reader's convenience. To this end, we use $\sym(\nabla u_i)=\C^{-1}\sigma_i+\sym p_i$ obtained from the elasticity relation and get
 $$\sym(\nabla (u_1-u_2))=\C^{-1}(\sigma_1-\sigma_2)+\sym(p_1-p_2)=0\,,$$
 and hence, from the first Korn's inequality (see e.g. \cite{NEFFKORN2002}), we get $\nabla(u_1-u_2)=0$ which implies that $u_1=u_2$. Therefore, we finally obtain
 $$u_1=u_2\,,\qquad \sigma_1=\sigma_2\,,\qquad(\underline{\gamma_1}=\underline{\gamma_2}\quad\Rightarrow\quad p_1=p_2)\,,\qquad \PL_{p_{\,1}}=\PL_{p_{\,2}},\qquad \underline{\eta_1}=\underline{\eta_2},$$ and thus the uniqueness
  of a weak/strong solution. \qed
  \begin{remark}\label{perfect-plasticity}{\rm It should be stressed that, in the proof, $k_2 > 0$ is necessary for uniqueness of the displacement field (and of the 
slip variables).
If $k_2 = 0$, we have perfect plasticity and multiple solutions involving displacement discontinuities along {\it slip lines}, as 
in conventional Hill's plasticity, are possible. The curl operator does not regularize such discontinuities since 
$\curl p$ may vanish in the presence of gradient of slip $\gamma^\alpha$ perpendicularly to the slip planes. 
It is shown in \cite{forestactamat98} that gradient models lead to finite width {\it kink bands} but still allow for 
{\it slip band} discontinuities, parallel to slip planes, in single crystals.}\end{remark}

\subsection{The model with linear kinematical hardening}\label{kinhard}
Here we consider the model where the isotropic hardening has been replaced with  linear kinematical hardening. 
\subsubsection{The description of the model}\label{description}
Here the free-energy density $\Psi$ is also given in the additively separated form as
\begin{eqnarray}\label{free-eng-kin}
\Psi(\nabla u,p,\PL_p,\Curl\PL_p):
&=&\underbrace{\Psi^{\mbox{\scriptsize lin}}_e(\bvarepsilon_e)}_{\mbox{\small elastic energy}}\qquad
+\qquad\underbrace{\Psi^{\mbox{\scriptsize lin}}_{\mbox{\scriptsize
 micro}}(p,\PL_p)}_{\mbox{\small micro
energy}}\\
&& \nonumber+\,\,\underbrace{\Psi^{\mbox{\scriptsize lin}}_{\mbox{\scriptsize
curl}}(\Curl \PL_p)}_{\mbox{\small defect-like 
energy (GND)}}\qquad+\hskip.4truecm\,\,\underbrace{\,\,\Psi^{\mbox{\scriptsize
lin}}_{\mbox{\scriptsize
 kin}}(\bvarepsilon_p)}_{\begin{array}{c}
\mbox{\small hardening energy (SSD)}\end{array}}\,,
 \end{eqnarray} where
 \begin{eqnarray}\label{free-eng-expr} \nonumber\Psi^{\mbox{\scriptsize
lin}}_e(\bvarepsilon_e) &:=&\frac12\,\la\bvarepsilon_e,\C_{\mbox{\tiny iso}} \bvarepsilon_e\ra=\frac12\,\la\sym(\nabla u-p),\C_{\mbox{\tiny iso}}\sym(\nabla u-p)\ra\,, \\
\Psi^{\mbox{\scriptsize lin}}_{\mbox{\tiny
 micro}}(p,\PL_p)&:=&\frac12\,\mu\,H_\chi\,\norm{p-\PL_p}^2\,,\qquad
 \Psi^{\mbox{\tiny lin}}_{\mbox{\tiny
curl}}(\Curl \PL_p)\,:=\,\frac12\,\mu\, L_c^2\norm{\Curl \PL_p}^2\,,\\ \nonumber\Psi^{\mbox{\scriptsize
lin}}_{\mbox{\tiny
 kin}}(\bvarepsilon_p)&:=&\frac12\,\mu\,k_1\norm{\bvarepsilon_p}^2=\frac12\,\mu\,k_1\norm{\sym p}^2
 \,.\end{eqnarray}
In this case, the equilibrium equation and the microcurl balance are obtained as in (\ref{EQL-crystal}) and in (\ref{EQL-micro-crystal}) respectively.\\
Now, the free-energy imbalance
$$\dot{\Psi}\leq\la\sigma,\nabla u\ra=\la\sigma,\dot{\bvarepsilon}_e\ra+\la\sigma,\dot{p}\ra$$ and the expansion of $\dot{\Psi}$ lead to the usual infinitesimal eleastic stress-strain 
relation
\begin{equation}
\nonumber \sigma =2\mu\, \sym(\nabla
u-p)+\lambda\, \tr(\nabla u-p)\id =2\mu\, (\sym(\nabla
u)-\bvarepsilon_p)+\lambda\, \tr(\nabla u)\id\label{elasticlaw-kin}
\end{equation}
and the local reduced dissipation inequality
\begin{equation}\label{diss-kinehard-crystal}
\la\Sigma_E,\dot{p}\ra\geq0\,,\end{equation}
where the non-symmetric Eshelby-type stress tensor in this case takes the form
\begin{equation}\label{eshelby-kin} \Sigma_E:=\sigma+\Sigma^{\mbox{\tiny lin}}_{\mbox{\tiny micro}} 
+\Sigma^{\mbox{\tiny lin}}_{\mbox{\tiny kin}}\end{equation}
with
\begin{eqnarray}
\Sigma^{\mbox{\tiny lin}}_{\mbox{\tiny micro}} &:=& \mu\,H_\chi(\PL_p-p)=-\,\mu\,L_c^2\Curl\Curl\PL_p,\label{sigma-micro}\\
\Sigma^{\mbox{\scriptsize lin}}_{\mbox{\scriptsize kin}} &:=&-\,\mu\,k_1\,\bvarepsilon_p=-\,\mu\,k_1\,\sym p\,.\label{sigma-kin}\end{eqnarray}
Two sources of kinematic hardening therefore arise in the model: the size--dependent contribution, $\Sigma^{\mbox{\tiny lin}}_{\mbox{\tiny micro}}$,
induced by strain gradient plasticity, and conventional size--independent linear kinematic hardening $\Sigma^{\mbox{\scriptsize lin}}_{\mbox{\scriptsize kin}}$.

Following the steps in the derivation of the strong formulation of the microcurl model with isotropic hardening in Section \ref{strong-iso-crystal}, we get the strong formulation in Table \ref{table:micro-kinehard-crystal} for the model with linear kinematical hardening.
\begin{table}[ht!]{\footnotesize\begin{center}
\begin{tabular}{|ll|}\hline &\qquad\\
 {\em Additive split of distortion:}& $\nabla u =e +p$,\quad $\bvarepsilon_e=\mbox{sym}\,e$,\quad $\bvarepsilon_p=\sym p$\\
{\em Equilibrium:} & $\mbox{Div}\,\sigma +f=0$ with
$\sigma=\C_{\mbox{\scriptsize iso}}\bvarepsilon^e=\C_{\mbox{\tiny iso}}(\sym\nabla u-\bvarepsilon_p)$\\
{\em Microbalance:} & $\mu\,L_c^2\,\Curl\Curl\PL_p=-\mu\,H_\chi\,(\PL_p-p)$,\\&\\
 {\em Free energy:} &
$\frac12\,\langle\C_{\mbox{\tiny iso}}\bvarepsilon_e,\bvarepsilon_e\rangle+\frac12\mu\,H_\chi\,\norm{p-\PL_p}^2$\\
&\qquad\quad$+\,\frac12\,\mu\,
L^2_c\,\norm{\Curl \PL_p}^2+\frac12\,\mu\, k_1\,\norm{\sym p}^2$\\&\\
{\em Yield condition:} &
$\phi(\Sigma_E):=\norm{\dev\Sigma_E}-\yieldzero\leq0$\\
 {\em where } & $\Sigma_E:=\sigma+\mu\,H_\chi\,(\PL_p-p)-\mu\,k_1\sym p$
 \\&\\{\em Dissipation inequality:} &
 $\langle\Sigma_E,\dot{p}\rangle\geq0$\\
 {\em Dissipation function:} &$\mathcal{D}_{\mbox{\scriptsize iso}}(q):=\yieldzero \norm{q}$\\
 {\em Flow rule in primal form:} &
 $\Sigma_E\in\partial \mathcal{D}_{\mbox{\scriptsize iso}}(\dot{p})$\\
{\em Flow rule in dual form:}
&$\dot{p}=\lambda\,\dsize\frac{\dev\Sigma_E}{\norm{\dev\Sigma_E}},\quad\qquad \lambda=\norm{\dot{p}}$\\&\\
{\em KKT conditions:} &$\lambda\geq0$, \quad $\phi(\Sigma_E,g)\leq0$,
\quad $\lambda\,\phi(\Sigma_E,g)=0$\\
 {\em Boundary conditions for $\PL_p$:} & $\PL_p\times{n}=0$ on
 $\Gamma_{\mbox{\scriptsize D}}$,\,\, $(\Curl \PL_p)\times{n}=0$ on $\partial\Omega\setminus\Gamma_{\mbox{\scriptsize D}}$\\
 {\em Function space for $\PL_p$:} & $\PL_p(t,\cdot)\in \mbox{H}(\mbox{Curl};\;\Omega,\,\BBR^{3\times 3})$\\
 \hline
\end{tabular}\caption{\footnotesize The microcurl model in single crystal gradient plasticity  with linear kinematical hardening. }\label{table:micro-kinehard-crystal}\end{center}}\end{table}

\subsubsection{The weak formulation of the model with linear kinematical hardening}\label{wk-kinehard-crystal}
The equilibrium and microbalance equations in weak form are
\begin{eqnarray}\label{weak-eq1-kin-crystal}
&&\hskip-1.5truecm\int_{\Omega}\la\C_{\mbox{\tiny iso}}(\sym\nabla
u-\bvarepsilon_p),\mbox{sym}(\nabla v-\nabla\dot{u})\ra\, dx=\int_\Omega
f(v-\dot{u})\,dx\,, \\
&&\label{EQL-microweak-kin-crystal}
\hskip-1.5truecm \int_\Omega\Bigl[\mu\,L_c^2\la\Curl\PL_p,\Curl X-\Curl\dot{\PL}_p\ra+\mu\,H_\chi\,\la\PL_p-p,X-\dot{\PL}_p\ra\Bigr]\,dx=0\,,
\end{eqnarray}
for every $v\in \SFV$ and $X\in \SFQ\,$  with $\SFV$ and  $\SFQ$ defined in (\ref{space-v}) and in (\ref{space-p}),  respectively.
\vskip.2truecm\noindent Now, the primal formulation of the flow rule ($\Sigma_E\in\partial\mathcal{D}_{\mbox{\scriptsize kin}}(\dot{p})$) in weak form reads  for every $q\in L^2(\Omega,\sL(3))$ as
\begin{eqnarray}\label{primal-wk-kin-crystal}\nonumber\int_\Omega\mathcal{D}_{\mbox{\tiny kin}}(q)dx-\int_\Omega\mathcal{D}_{\mbox{\tiny kin}}(\dot{p})dx&\geq&\int_\Omega\la\Sigma_E,q-\dot{p}\ra\,dx\\
&=&\int_\Omega\la\C_{\mbox{\tiny iso}}\sym(\nabla u-p),\sym(q-\dot{p}\ra\,dx\\
&& \nonumber+\int_\Omega\Bigl[\la\mu\,H_\chi\,(\PL_p-p)-\mu\,k_1\sym p,q-\dot{p}\ra\Bigr]dx\,.\end{eqnarray}

Now adding up (\ref{weak-eq1-kin-crystal}), (\ref{EQL-microweak-kin-crystal})  and (\ref{primal-wk-kin-crystal}) we
get the following weak formulation of the microcurl model of single crystal strain gradient plasticity with linear kinematical hardening in the form of a variational inequality:
\begin{eqnarray}\label{weak-form-iso}
\nonumber&&\hskip-.8truecm\int_\Omega\Bigl[\la\C_{\mbox{\tiny iso}}\sym(\nabla
u-p),\sym(\nabla v-q)-\sym(\nabla\dot{u}-\dot{p})\ra+ \mu\,L_c^2\la\Curl\PL_p,\Curl X-\Curl\dot{\PL}_p\ra\\
\nonumber&&\qquad+\,\,\mu\,H_\chi\,\la\PL_p-p, (X-q)-(\dot{\PL}_p-\dot{p})\ra+\mu\,k_1\,\la\sym p,\sym(q-\dot{p})\ra\Bigr]dx\\
&&\qquad  +
\int_\Omega\mathcal{D}_{\mbox{\scriptsize kin}}(q)\,dx
-\int_\Omega\mathcal{D}_{\mbox{\scriptsize kin}}(\dot{p})\,dx\,\geq\, \int_\Omega
f\,(v-\dot{u})\,dx\,.\end{eqnarray}
That is, setting $\SFZ:=\SFV\times \SFP\times\SFQ$ with $\SFV$, $\SFP$ and $\SFQ$ defined in (\ref{space-v})-(\ref{space-p}) and their norms in (\ref{norm-Z}), we get the problem of the form: Find $w=(u,p,\PL_p)\in  \SFH^1(0,T;Z)$
such that $w(0)=0$ and 
\begin{equation}\label{wf-kin-crystal}
\ba(\dot{w},z-w)+j(z)-j(\dot{w})\geq \langle
\ell,z-\dot{w}\rangle\mbox{ for every } z\in \SFZ\mbox{ and for a.e.
}t\in[0,T]\,,\end{equation}
where
\begin{eqnarray}
\ba(w,z)&:=&\int_\Omega\Bigl[\la\C_{\mbox{\tiny iso}}\sym(\nabla
u-p),\sym(\nabla v-q)\ra+ \mu\,L_c^2\la\Curl\PL_p,\Curl X\ra\\
\nonumber&&\hskip2truecm+\,\,\mu\,H_\chi\,\la\PL_p-p, X-q\ra+\mu\,k_1\,\la\sym p,\sym q\ra\,\Bigr]dx\,,
\label{bilin-ikin-crystal}
\\\nonumber\\
j(z)&:=&\int_\Omega\mathcal{D}_{\mbox{\tiny kin}}(q)\,dx\,,\label{functional-kin-crystal}\\
\langle \ell,z\rangle&:=&\int_\Omega
f\,v\,dx\,,\label{lin-form-kin-crystal}\end{eqnarray} for
$w=(u,p,\PL_p)$ and $z=(v,q,\QL)$ in
$\SFZ$.
\subsubsection{Existence and uniqueness for the model with linear kinematic hardening}\label{exist-linkin-cryst}
In order to show the existence and uniqueness  for the problem in  (\ref{wf-kin-crystal})-(\ref{lin-form-kin-crystal}) using  \cite[Theorem 6.15]{Han-ReddyBook}, we only need to show here that  the bilinear form $\ba$ is $\SFZ$-coercive. However, this is obtained following a different approach. We will make use of the following result.

\begin{lemma}\label{new-norm-crystal}
The mapping $\norm{\cdot}_*:\SFP\times\SFQ\to[0,\infty)$ defined by
\begin{equation}\label{new-norm1} \norm{(q,X)}_*^2:=\norm{q-X}^2_{L^2}+\norm{\sym q}^2_{L^2}+\norm{\Curl X}^2_{L^2}\end{equation} is  a norm on $\SFP\times \SFQ$ equivalent to the norm defined by \begin{equation}\label{usual-norm}\norm{(q,X)}^2_{P\times Q}=\norm{q}^2_{L^2}+\norm{X}^2_{\mbox{\scriptsize $H$}(\mbox{\scriptsize
$\Curl$};\Omega)}\,.\end{equation}
\end{lemma}
{\bf Proof.} To show that $\norm{\cdot}_*$ is a norm on $\SFP\times\SFQ$, we only  check the vanishing property of a norm since the other properties are trivially satisfied. The vanishing property is obtained through the Korn-type inequality for incompatible tensor fields established  in \cite{NPW2011-1, NPW2012-1, NPW2012-2,NPW2014}, namely
\begin{equation}\label{korn-incomp-crystal}
\norm{X}^2_{L^2}\leq C\,(\norm{\sym X}^2_{L^2}+\norm{\Curl X}^2_{L^2})\qquad\forall X\in\SFQ:=\mbox{H}_0(\mbox{Curl};\,\Omega,\,\Gamma_{\mbox{\scriptsize D}},\BBR^{3\times 3})\,.\end{equation}
In fact, let $(q,X)\in\SFP\times\SFQ$ be such that $\norm{(q,X)}_*=0$, that is, $q=X$, $\sym q=0$ and $\Curl  X=0$. Thus we get $\sym X=0$ and $\Curl X=0$. From (\ref{korn-incomp-crystal}), we then get $X=0$ and hence also $q=0$.

Now to show that both norms are equivalent, we will  first show that $(\SFP\times\SFQ,\norm{\cdot}_*)$ is a Banach space. To this aim, let $(q_n,X_n)$ be a Cauchy sequence in $(\SFP\times\SFQ,\norm{\cdot}_*)$. Hence, the sequences $(q_n-X_n)$, $(\sym q_n)$ and $(\Curl X_n)$ are all Cauchy in $L^2(\Omega,\mathbb{R}^{3\times 3})$ and therefore, there exist $A,\,B,\,C\in L^2(\Omega,\mathbb{R}^{3\times 3})$ such that 
\begin{equation}\label{limit-cauchy-seq}q_n-X_n\to A,\qquad \sym q_n\to B,\qquad\mbox{and}\qquad\Curl X_n\to C\,.\end{equation}
Thus, $$\sym X_n=\sym q_n-\sym(q_n-X_n)\to B-\sym A=\sym(B-A)\mbox{\, and \, }\Curl X_n\to C\,.$$ Hence, it follows from the inequality (\ref{korn-incomp-crystal}) that $(X_n)$ is a Cauchy sequence in $(\SFQ,\norm{\cdot}_{\mbox{\scriptsize H}(\mbox{\scriptsize
Curl};\Omega)})$. Hence, there exists $X\in\SFQ$ such that 
$$X_n\to X\quad\mbox{ and }\quad  \Curl X_n\to\Curl X=C\mbox{ in }L^2(\Omega,\mathbb{R}^{3\times 3})\,.$$
Now $q_n=X_n+(q_n-X_n)\to X+A$ and $\sym q_n\to\sym(X+A)=B$. Therefore,
\begin{eqnarray*}\norm{(q_n,X_n)-(X+A,X)}_*^2&=&\norm{(q_n-(X+A),X_n-X)}_*^2\,=\,\norm{(q_n-X_n)-A}^2_{L^2}\\
&&\quad +\norm{\sym q_n-\sym(X+A)}^2_{L^2} +\norm{\Curl X_n-\Curl X}^2_{L^2}\to0\,.\end{eqnarray*}
So, the sequence $(q_n,X_n)$ converges to $(X+A,X)$ in $(\SFP\times\SFQ,\norm{\cdot}_*)$.\\ Since the two normed spaces   $(\SFP\times\SFQ,\norm{\cdot}_{P\times Q})$ and  $(\SFP\times\SFQ,\norm{\cdot}_*)$
 are Banach and the identity mapping 
 $$\mbox{Id}:(\SFP\times\SFQ,\norm{\cdot}_{P\times Q})\to(\SFP\times\SFQ,\norm{\cdot}_*)$$ is linear and continuous, then as a consequence of the open mapping theorem, we find that
 $$\mbox{Id}:(\SFP\times\SFQ,\norm{\cdot}_*)\to(\SFP\times\SFQ,\norm{\cdot}_{P\times Q})$$ is also linear and continuous. Therefore, the  two norms  $\norm{\cdot}_*$ and $\norm{\cdot}_{P\times Q}$ are equivalent and this completes the proof of the lemma.\hfill\qed
 \vskip.2truecm\noindent We are now in a position to prove that the bilinear form $\ba$ in (\ref{bilin-ikin-crystal}) is $\SFZ$-coercive.
 \begin{lemma}\label{elliptic-a}
 There exists a positive constant $C$ such that $\ba(z,z)\geq C\norm{z}^2_Z$ for every $z\in \SFZ$.\end{lemma}
 {\bf Proof.} Let $z=(v,q,X)\in\SFZ=\SFV\times\SFP\times\SFQ$

\begin{eqnarray}
\nonumber\ba(z,z)& \geq&  m_0\norm{\sym(\nabla v-q)}^2_{L^2}\mbox{
(from (\ref{ellipticityC}))} +\, \mu\,H_\chi\,\norm{X-q}^2_{L^2}\\
\nonumber&&\qquad\qquad\hskip4.1truecm+\,\mu\,
L_c^2\norm{\Curl X}_{L^2}^2 +\mu\, k_1\,\norm{\sym q}_{L^2}^2 \\
\nonumber&=&m_0\bigl[\norm{\sym\nabla v}^2_{L^2} +\norm{\sym q}^2_2
-2\la\sym\nabla v, \sym q\ra\bigr]+\mu\,H_\chi\,\norm{X-q}^2_{L^2}\\
\nonumber&&\qquad\qquad\hskip4.1truecm  +\,\mu L_c^2\,\norm{\Curl X}_{L^2}^2
+\mu\,k_1\norm{\sym q}_{L^2}^2\\
\nonumber&\geq&m_0\,(1-\theta)\norm{\sym\nabla v}^2_{L^2} +\Bigl[m_0\Bigl(1-\frac1\theta\Bigr)+\mu\,k_1\Bigr]\norm{\sym q}^2_2+\mu\,H_\chi\,\norm{X-q}^2_{L^2}\\
\nonumber&&\qquad\qquad\hskip4.1truecm  +\,\mu L_c^2\,\norm{\Curl X}_{L^2}^2\,.
\end{eqnarray}
Now,  since the hardening constant $k_1>0$, we choose $\theta$ such that 
$$\frac{m_0}{m_0+\mu\,k_1}<\theta<1\,,$$ and  using Korn's first inequality (see e.g. \cite{NEFFKORN2002}) and Lemma \ref{new-norm-crystal}, we then get two constants $C=C(\theta,m_0,\mu,H_\chi,k_1,L_c,\Omega)>0$ and $C'=C'(\theta,m_0,\mu,H_\chi,k_1,L_c,\Omega)>0$ such that 
$$\ba(z,z) \geq C'\Bigl[\norm{\nabla v}^2_{L^2}+\norm{(q,X)}^2_*\Bigr]\geq C\Bigl[\norm{\nabla v}^2_{L^2}+\norm{q}^2_{L^2}+\norm{X}^2_{\mbox{\scriptsize $H$}(\mbox{\scriptsize
$\Curl$};\Omega)}\Bigr]=C\norm{z}^2_Z\,.\qquad\mbox{\qed}$$

\section{The microcurl model in polycrystalline gradient plasticity}\label{poly}
\subsection{The case with isotropic hardening}\label{micro-poly-iso}

The  free-energy density $\Psi$ is given in the additively separated form
\begin{eqnarray}\label{freeng-poly-iso}
\Psi(\nabla u,\bvarepsilon_p,\PL_p,\Curl\PL_p,\eta_p):
&=&\underbrace{\Psi^{\mbox{\scriptsize lin}}_e(\bvarepsilon_e)}_{\mbox{\small elastic energy}}\qquad
+\qquad\underbrace{\Psi^{\mbox{\scriptsize lin}}_{\mbox{\scriptsize
 micro}}(\bvarepsilon_p,\PL_p)}_{\mbox{\small micro
energy}}\\
&& \nonumber+\,\,\underbrace{\Psi^{\mbox{\scriptsize lin}}_{\mbox{\scriptsize
curl}}(\Curl \PL_p)}_{\mbox{\small defect-like
energy (GND)}}\qquad+\hskip.4truecm\,\,\underbrace{\,\,\Psi^{\mbox{\scriptsize
lin}}_{\mbox{\scriptsize
 iso}}(\eta_p)}_{\begin{array}{c}
\mbox{\small hardening energy (SSD)}\end{array}}\,,
 \end{eqnarray} where
 \begin{eqnarray}\label{freeng-expr-poly-iso} \nonumber\Psi^{\mbox{\scriptsize
lin}}_e(\bvarepsilon_e) &:=&\frac12\,\la\bvarepsilon_e,\C_{\mbox{\scriptsize iso}} \bvarepsilon_e\ra=\frac12\,\la\sym\nabla u-\bvarepsilon_p,\C_{\mbox{\scriptsize iso}}(\sym\nabla u-\bvarepsilon_p)\ra\,, \\
\Psi^{\mbox{\scriptsize lin}}_{\mbox{\scriptsize
 micro}}(\bvarepsilon_p,\PL_p)&:=&\frac12\,\mu\,H_\chi\,\norm{\bvarepsilon_p-\sym\PL_p}^2\,,\quad
 \Psi^{\mbox{\scriptsize lin}}_{\mbox{\scriptsize
curl}}(\Curl \PL_p)\,:=\,\frac12\,\mu\, L_c^2\norm{\Curl \PL_p}^2\,,\\ \nonumber\Psi^{\mbox{\scriptsize
lin}}_{\mbox{\scriptsize
 iso}}(\eta_p)&:=&\frac12\mu\,k_2|\eta_p|^2
 \,.\end{eqnarray}
Here, $\eta_p$  is the isotropic hardening variable. \\
It should be noted that there is no constraint on the 
skew--symmetric part of the microdeformation, $\skew \PL_p$, in (\ref{freeng-expr-poly-iso}), due to the fact that
no plastic spin is considered in the original plasticity model for $\skew \PL_p$ to be compared with. It will be shown
that, in spite of that, no indeterminacy of $\skew \PL_p$ arises in the formulation\footnote{No simple characterization of $\skew\PL_p$ for $H_\chi\to\infty$ is known at present.}. This represents
the most straightforward microcurl extension of a phenomenological polycrystal plasticity model. 
\subsubsection{The balance equations.}\label{Equi-poly}
As in Section \ref{Equi}, we have the balance equations:
\begin{eqnarray}
\label{EQL-macro}
\div \sigma + f &=& 0 \,\qquad\hskip5truecm\mbox{(macroscopic balance)}\,,\\
\label{EQL-micro}
\mu\,L^2_c\Curl\Curl\PL_p&=&-\,\mu\,H_\chi\,(\sym\PL_p-\bvarepsilon_p)\in\mbox{Sym}(3)\,\qquad\mbox{(microbalance)}\,,
\end{eqnarray}
 where (\ref{EQL-micro}) is supplemented by the boundary conditions
\begin{equation}\label{bc-poly}
\PL_p\times n|_{\Gamma_{\mbox{\scriptsize D}}}=0\qquad\mbox{ and }\qquad (\Curl\PL_p)\times n|_{\partial\Omega\setminus\Gamma_{\mbox{\scriptsize D}}}=0\,.
\end{equation}

\subsubsection{The derivation of the dissipation inequality.}\label{deriv-diss-ineq-poly}

The local free-energy imbalance states that
\begin{equation}
\dot{\Psi} - \la\sigma,\dot{\bvarepsilon}_e\ra - \la\sigma,\dot{\bvarepsilon}_p\ra  \leq 0\
. \label{2ndlaw-poly}
\end{equation}
Now we expand the first term, substitute (\ref{freeng-poly-iso}) and get
\begin{equation}\label{2ndlaw-poly-iso}
\la\C_{\mbox{\scriptsize iso}}\,\bvarepsilon_e-\sigma,\dot{\bvarepsilon}_e\ra-\la\sigma,\dot{\bvarepsilon}_p\ra-\mu\,H_\chi\,\la\sym\PL_p-\bvarepsilon_p,\dot{\bvarepsilon}_p\ra 
+\mu\,k_2\,\eta_p\,\dot{\eta}_p\leq0\,.
\end{equation}
Since the inequality (\ref{2ndlaw-poly-iso}) must be satisfied for whatever elastic-plastic deformation mechanism, inlcuding purely elastic ones (for which $\dot{\eta}_p=0$, $\dot{\varepsilon}_p=0)$, then it implies the infinitesimal stress-strain 
relation
\begin{equation}
\sigma = \C_{\mbox{\scriptsize iso}}\,\varepsilon_e=2\mu\, (\sym\nabla
u-\varepsilon_p)+\lambda\, \tr(\sym\nabla u-\varepsilon_p)\id  \label{elasticlaw-poly}
\end{equation}
and the local reduced dissipation inequality
\begin{equation}\label{diss-poly-iso1}
-\la\sigma,\dot{\bvarepsilon}_p\ra-\mu\,H_\chi\la\sym\PL_p-\bvarepsilon_p,\dot{\bvarepsilon}_p\ra+\mu\,k_2\,\eta_p\,\dot{\eta}_p\leq0\,.\end{equation}
That is,
\begin{equation}\label{diss-poly-iso2}
\la\sigma+\mu\,H_\chi\,(\sym\PL_p-\bvarepsilon_p),\dot{\bvarepsilon}_p\ra -\mu\,k_2\,\eta_p\,\dot{\eta}_p
\geq0\,,\end{equation}
which can also be written in compact form as 
\begin{equation}\label{diss-poly-iso3}
\la\Sigma_p,\dot{\Gamma}_p\ra\geq0\end{equation} where
\begin{equation}\label{gen-var-poly-iso}\Sigma_p:=(\Sigma_E,g)\qquad\mbox{and}\qquad \Gamma_p=(\bvarepsilon_p,\eta_p)
\end{equation}
with $\Sigma_E$ being a symmetric Eshelby-type stress tensor and $g$ being a thermodynamic force-type variable conjugate to $\dot{\eta}_p$ and defined as
\begin{eqnarray} \Sigma_E&:=&\sigma+\mu\,H_\chi\,(\sym\PL_p-\bvarepsilon_p)\,\,=\,\,\sigma-\mu\,L_c^2\,\Curl\Curl\PL_p,\\
 g&:=&-\,\mu\,k_2\,\eta_p.\end{eqnarray}

\subsubsection{The flow rule}\label{flow-law}
We consider a yield function
defined by
\begin{equation}
\label{yield-funct-kinematic} \phi(\Sigma_p):= \norm{\dev\Sigma_{\mbox{\scriptsize E}}}+g - \yieldzero\quad\mbox{ for }\Sigma_p=(\Sigma_E,g)\,.
\end{equation}
So the set
of admissible (elastic) generalized stresses is defined as
\begin{equation}\label{admiss-stress-kin}
\mathcal{K}:=\left\{\Sigma_p=(\Sigma_{\mbox{\scriptsize
E}},g)\,\,|\,\,\phi(\Sigma_p)
   \leq0,\,\,g\leq0\right\}\,.\end{equation}
The  principle of maximum dissipation gives the normality
law\begin{equation}\label{normalcone} \dot{\Gamma}_p\in
N_{\mathcal{K}}(\Sigma_p)\,,\end{equation}
where $\dsize N_{\mathcal{K}}(\Sigma_p)$
denotes the normal cone to $\mathcal{K}$ at
$\dsize\Sigma_p$, which is the set of
generalised strain rates $\dot{\Gamma}_p$ that satisfy
\begin{equation}
\la\overline{\Sigma} - \Sigma_p,\dot{\Gamma}_p\ra \leq
0\ \quad \mbox{for all}\ \overline{\Sigma} \in {\mathcal{K}}\ .
\label{normality2}
\end{equation} Notice that $N_{\mathcal{K}}=\partial\Chi_{\mathcal{K}}$ where
$\Chi_{\mathcal{K}}$ denotes the indicator function of the set $\mathcal{K}$
and $\partial\Chi_{\mathcal{K}}$ denotes the subdifferential of the function
$\Chi_{\mathcal{K}}$.\\ Whenever the yield surface $\partial\mathcal{K}$ is
smooth at $\dsize\Sigma_p$ then
$$\dot{\Gamma}_p\in
N_{\mathcal{K}}(\Sigma_p)\quad\Rightarrow\quad\exists\lambda\mbox{ such that }
\dot{\bvarepsilon}_p=\lambda\,\frac{\dev\Sigma_{\mbox{\scriptsize E}}}{\norm{\dev\Sigma_{\mbox{\scriptsize E}}}}\quad\mbox{and}\quad\dot{\eta}_p=\lambda=\norm{\dot{\bvarepsilon}_p}$$ with the Karush-Kuhn
Tucker conditions: $\lambda\geq0$, $\phi(\Sigma_p)\leq0$ and $\lambda\,\phi(\Sigma_p)=0$\,.\\
Using convex analysis (Legendre-transformation) we find that
\begin{equation}\label{primalflowlaw}\underbrace{\dot{\Gamma}_p\in
\partial\Chi_{\mathcal{K}}(\Sigma_p)}_{\mbox{\bf flow rule in its dual formulation}}\quad\Leftrightarrow\quad \underbrace{\Sigma_p\in
\partial \Chi^*_{\mathcal{K}} (\dot{\Gamma}_p)\,}_{\mbox{\bf flow rule in its primal formulation}}
\end{equation} where $\Chi^*_{\mathcal{K}}$ is the Fenchel-Legendre dual of the function $\Chi_{\mathcal{K}}$ denoted in this context by $\mathcal{D}_{\mbox{\scriptsize iso}}$,
 the one-homogeneous dissipation function for rate-independent processes. That is, for every $\Gamma=(q,\beta)$,
\begin{eqnarray}\label{dissp-function-iso}
\nonumber\mathcal{D}_{\mbox{\scriptsize iso}}(\Gamma)&=&\sup\graffe{\la\Sigma_p,\Gamma\ra\,\,|\,\,\Sigma_p\in\mathcal{K}}\\
\nonumber&=&\sup \graffe{\la\Sigma_E,q\ra +g\beta \,\,|\,\
\phi(\Sigma_E,g)\leq0,\,\,g\leq0}\\&=&\left\{\begin{array}{ll}
      \yieldzero\,\norm{q} &\mbox{ if } \norm{q}\leq\beta\,,\\
      \infty &\mbox{ otherwise.}\end{array}\right.
\end{eqnarray} We get from the definition of the subdifferential ($\Sigma_p \in
\partial \Chi^*_{\mathcal{K}} (\dot{\eta}_p)$) that,
\begin{equation}
\mathcal{D}_{\mbox{\scriptsize iso}} (\Gamma) \geq \mathcal{D}_{\mbox{\scriptsize iso}}(\dot{\Gamma}_p) +
\la\Sigma_p,\Gamma-\dot{\Gamma}_p\ra\quad \mbox{for any
} \Gamma.\label{dissineq}\end{equation}
That is,
\begin{equation} \mathcal{D}_{\mbox{\scriptsize iso}}(q,\beta)\geq \mathcal{D}_{\mbox{\scriptsize iso}}(\dot{\varepsilon}_p,\dot{\eta}_p)+\la\Sigma_E,q-\dot{\varepsilon}_p\ra+g(\beta-\dot{\eta}_p)\quad \mbox{for any
} (q,\beta).\label{dissinequality2}\end{equation}
\vskip.3truecm\noindent In the next sections, we present as in the case of signle-crystal gradient  plasticity, a complete mathematical analysis of the model including  both strong and weak fomrulations as well as a corresponding
existence result.
\subsubsection{Strong formulation
of the model}\label{strong-iso-poly} To summarize, we have obtained the following strong
formulation for the microcurl model in the poycrystalline infinitesimal gradient plasticity setting with
 isotropic  hardening. Given $f\in \SFH^1(0,T;L^2(\Omega,\mathbb{R}^3))$, the goal is to find:
\begin{itemize}\item[(i)] the displacement $u\in \SFH^1(0,T; H^1_0(\Omega,{\Gamma_{\mbox{\scriptsize D}}},\mathbb{R}^3))$,
\item[(ii)] the infinitesimal plastic strain $\bvarepsilon_p\in
\SFH^1(0,T;L^2(\Omega, \mbox{Sym}(3)\cap\sL(3)))$, the infinitesimal micro-distortion $\PL_p$ with  $\sym\PL_p\in\SFH^1(0,T;L^2(\Omega, \mbox{Sym}(3)\cap\sL(3)))$, $\Curl\PL_p\in \SFH^1(0,T;
L^2(\Omega,\BBR^{3\times 3}))$ and $ \Curl\Curl\PL_p \in \SFH^1(0,T;
L^2(\Omega,\BBR^{3\times 3}))$\end{itemize}
 such that the content of Table \ref{table:micro-isohard-poly} holds.
\begin{table}[ht!]{\footnotesize\begin{center}
\begin{tabular}{|ll|}\hline &\qquad\\
 {\em Additive split of strain:}&  $\nabla u =e +p$,\quad $\bvarepsilon_e=\mbox{sym}\,e$,\quad $\bvarepsilon_p=\sym p$\\
{\em Equilibrium:} & $\mbox{Div}\,\sigma +f=0$ with
$\sigma=\C_{\mbox{\scriptsize iso}}\bvarepsilon^e=\C_{\mbox{\scriptsize iso}}(\sym\nabla u-\bvarepsilon_p)$\\
{\em Microbalance:} & $\mu\,L_c^2\,\Curl\Curl\PL_p=-\mu\,H_\chi\,(\sym\PL_p-\bvarepsilon_p)$,\\&\\
 {\em Free energy:} &
$\frac12\,\langle\C.\bvarepsilon^e,\bvarepsilon^e\rangle+\frac12\mu\,H_\chi\,\norm{\bvarepsilon_p-\sym\PL_p}^2$\\
&\qquad\qquad $+\frac12\,\mu\,
L^2_c\,\norm{\Curl \PL_p}^2+\frac12\,\mu\, k_2\,|\eta_p|^2$\\&\\
{\em Yield condition:} &
$\phi(\Sigma_E):=\norm{\dev\Sigma_E}+g-\yieldzero\leq0$\\
 {\em where } & $\Sigma_E:=\sigma+\mu\,H_\chi\,(\sym\PL_p-\bvarepsilon_p)$,\,\,\quad $g=-\mu\,k_2\,\eta_p$
 \\&\\{\em Dissipation inequality:} &
 $\langle\Sigma_E,\dot{\bvarepsilon}_p\rangle +g\dot{\eta}_p\geq0$\\
 {\em Dissipation function:} &$\mathcal{D}_{\mbox{\scriptsize iso}}(q,\beta):=\left\{\begin{array}{ll}\yieldzero \norm{q} &\mbox{ if }\norm{q}\leq\beta,\\ \infty &\mbox{ otherwise}\end{array}\right.$\\
 {\em Flow rule in primal form:} &
 $(\Sigma_E,g)\in\partial \mathcal{D}_{\mbox{\scriptsize iso}}(\dot{\bvarepsilon}_p,\dot{\eta}_p)$\\&\\
{\em Flow rule in dual form:}
&$\dot{\bvarepsilon}_p=\lambda\,\dsize\frac{\dev\Sigma_E}{\norm{\dev\Sigma_E}},\quad\qquad \dot{\eta}_p=\lambda=\norm{\dot{\bvarepsilon}_p}$\\&\\
{\em KKT conditions:} &$\lambda\geq0$, \quad $\phi(\Sigma_E,g)\leq0$,
\quad $\lambda\,\phi(\Sigma_E,g)=0$\\
 {\em Boundary conditions for $\PL_p$:} & $\PL_p\times{n}=0$ on
 $\Gamma_{\mbox{\scriptsize D}}$,\,\, $(\Curl \PL_p)\times{n}=0$ on $\partial\Omega\setminus\Gamma_{\mbox{\scriptsize D}}$\\
 {\em Function space for $\PL_p$:} & $\PL_p(t,\cdot)\in \mbox{H}(\mbox{Curl};\;\Omega,\,\BBR^{3\times 3})$\\
 \hline
\end{tabular}\caption{\footnotesize The microcurl model in polycrystalline gradient plasticity  with isotropic hardening.  The boundary condition on $\PL_p$ necessitates at least $\PL_p\in
\mbox{H}(\mbox{Curl};\,\Omega,\,\BBR^{3\times 3})$. This is proven to be the case in the next sections through a weak formulation of the model as a variational inequality.}\label{table:micro-isohard-poly}\end{center}}\end{table}
\subsubsection{Weak formulation of the model}\label{wf-iso-poly} Assume that the problem in Section \ref{strong-iso-poly} has a solution

$(u,\bvarepsilon_p,\PL_p,\eta_p)$.  Let $v\in H^1(\Omega,\mathbb{R}^3)$
with $v_{|\Gamma_D}=0$. Multiply the equilibrium equation with
$v-\dot{u}$ and integrate in space by parts and use the
symmetry of $\sigma$ and the elasticity relation to get
\begin{equation}\label{weak-eq1}
\int_{\Omega}\la\C_{\mbox{\scriptsize iso}}(\sym\nabla
u-\bvarepsilon_p)),\mbox{sym}\nabla v-\sym\nabla\dot{u}\ra\, dx=\int_\Omega
f(v-\dot{u})\,dx\, .
\end{equation}
Now,
for any $X\in C^\infty(\overline{\Omega},\sL(3))$ such that
$X\times\,n=0$ on $\Gamma_{\mbox{\scriptsize D}}$ 
we integrate (\ref{EQL-micro}) over $\Omega$, integrate by
parts the term with Curl\,Curl using the boundary conditions
$$(X-\dot{\PL}_p)\times\,n=0\mbox{ on }\Gamma_{\mbox{\scriptsize D}},\qquad
\mbox{Curl}\,\PL_p\times\,n=0\mbox{ on }
\partial\Omega\setminus\Gamma_{\mbox{\scriptsize D}}$$ and get
\begin{equation}\label{EQL-microweak}
\int_\Omega\Bigl[\mu\,L_c^2\la\Curl\PL_p,\Curl X-\Curl\dot{\PL}_p\ra-\mu\,H_\chi\,\la\bvarepsilon_p-\sym\PL_p,\sym X-\sym\dot{\PL}_p\ra\Bigr]\,dx=0\,.
\end{equation}
Moreover,
for any $q\in C^\infty(\overline{\Omega},\sL(3))$  and any $\beta\in L^2(\Omega)$,
we integrate (\ref{dissinequality2}) over $\Omega$ and get

\begin{eqnarray}\label{dissinequality2-wk}
\nonumber &&\hskip-1truecm \int_\Omega\mathcal{D}_{\mbox{\scriptsize iso}}(q,\beta)\,dx
-\int_\Omega\mathcal{D}_{\mbox{\scriptsize iso}}(\dot{\bvarepsilon}_p,\dot{\eta}_p)\,dx -\int_\Omega
\la\C_{\mbox{\scriptsize iso}}(\sym(\nabla u)-\bvarepsilon_p),q-\dot{\bvarepsilon}_p\ra\,dx\\
&&\quad+\int_\Omega\Bigl[\mu\,H_\chi\,\la\bvarepsilon_p-\sym\PL_p,q-\dot{\bvarepsilon}_p\ra
+\mu\,\alpha_1\,\eta_p\,(\beta-\dot{\eta}_p)\,\Bigr]dx\geq0\,.\end{eqnarray}

Now adding up (\ref{weak-eq1}), (\ref{EQL-microweak})  and (\ref{dissinequality2-wk}) we
get the following weak formulation of the problem in Section
\ref{strong-iso-poly} in the form of a variational inequality:
\begin{eqnarray}\label{weak-form-iso-poly}
\nonumber&&\hskip-.8truecm\int_\Omega\Bigl[\la\C_{\mbox{\scriptsize iso}}(\sym\nabla
u-\bvarepsilon_p),(\sym\nabla v-q)-(\sym\nabla\dot{u}-\dot{\bvarepsilon}_p)\ra+ \mu\,L_c^2\la\Curl\PL_p,\Curl X-\Curl\dot{\PL}_p\ra\\
&&\qquad+\,\,\mu\,H_\chi\,\la\sym\PL_p-\bvarepsilon_p, (\sym X-q)-(\sym\dot{\PL}_p-\dot{\bvarepsilon}_p)\ra+\mu\,k_2\,\eta_p\,(\beta-\dot{\eta}_p)\,\Bigr]dx\\
\nonumber &&\qquad  +
\int_\Omega\mathcal{D}_{\mbox{\scriptsize iso}}(q,\beta)\,dx
-\int_\Omega\mathcal{D}_{\mbox{\scriptsize iso}}(\dot{\bvarepsilon}_p,\dot{\eta}_p)\,dx\,\geq\, \int_\Omega
f\,(v-\dot{u})\,dx\,\end{eqnarray}

\subsubsection{Existence result for the weak formulation}\label{exist-poly-iso}
As in the case of single-crystal gradient plasticity in Section \ref{exist-crystal-iso}, the existence result for the weak formulation (\ref{weak-form-iso-poly}) is obtaned through the abstract machinery developed in \cite{Han-ReddyBook} for mathematical problems in geometrically linear
classical plasticity. To this aim,
(\ref{weak-form-iso-poly}) is written as the variational inequality of
the second kind: find $w=(u,\bvarepsilon_p,\PL_p,\eta_p)\in \SFH^1(0,T;Z)$
such that $w(0)=0$,  $\dot{w}(t)\in W$ for a.e. $t\in[0,T]$ and
\begin{equation}\label{wf}
\ba(\dot{w},z-w)+j(z)-j(\dot{w})\geq \langle
\ell,z-\dot{w}\rangle\mbox{ for every } z\in W\mbox{ and for a.e.
}t\in[0,T]\,,\end{equation} where $Z$ is a suitable Hilbert space
and $W$ is some closed, convex subset of $Z$ to be constructed
later,
\begin{eqnarray}
\ba(w,z)&=&\int_\Omega\Bigl[\la\C_{\mbox{\scriptsize iso}}(\sym\nabla
u-\bvarepsilon_p),\sym\nabla v-q\ra+ \mu\,L_c^2\la\Curl\PL_p,\Curl X\ra\\
\nonumber&&\hskip2truecm+\,\,\mu\,H_\chi\,\la\sym\PL_p-\bvarepsilon_p, \sym X-q\ra+\mu\,k_2\,\eta_p\,\beta\,\Bigr]dx\,,
\label{bilin-iso-spin}
\\\nonumber\\
j(z)&=&\int_\Omega\mathcal{D}_{\mbox{\scriptsize iso}}(q,\beta)\,dx\,,\label{functional-isospin}\\
\langle \ell,z\rangle&=&\int_\Omega
f\,v\,dx\,,\label{lin-form}\end{eqnarray} for
$w=(u,\bvarepsilon_p,\PL_p,\eta_p)$ and $z=(v,q,X,\beta)$ in
$\SFZ$.
\vskip.2truecm\noindent
The Hilbert space $\SFZ$ and the closed convex subset $\SFW$ are
constructed in such a way that the functionals $\ba$, $j$ and
$\ell$ satisfy the assumptions in the abstract result in
\cite[Theorem 6.19]{Han-ReddyBook}. The key issue here is the
coercivity of the bilinear form $\ba$ on the set $\SFW$, that is, $a(z,z)\geq C\norm{z}^2_Z$ for every $z\in\SFZ$ and for some $C>0$.
\vskip.2truecm\noindent
We let
\begin{eqnarray}
\SFV&:=&\mathsf{H}^1_0(\Omega,{\Gamma_{\mbox{\scriptsize
D}}},\mathbb{R}^3)=\{v\in \mathsf{H}^1(\Omega,\mathbb{R}^3)\,|\,
v_{|\Gamma_{\mbox{\scriptsize
D}}}=0\}\,,\label{space-v}\\
 \SFP&:=&\mbox{L}^2(\Omega,\,\sL(3)\cap\mbox{Sym}(3))\,,\label{space-epsp}\\
\SFQ&:=&\mbox{H}_0(\mbox{Curl};\,\Omega,\,\Gamma,\sL(3))\,,\label{space-p}\\
\Lambda&:=& L^2(\Omega)\,,\label{space-hardvar}\\
   \SFZ&:=&\SFV\times \SFP\times\SFQ\times\Lambda\,,\label{product-space}\\
   \SFW&:=&\bigl\{z=(v,q,X,\beta)\in\SFZ\,\,|\,\,\norm{q}\leq\beta\bigr\}\,,\label{set-W}
\end{eqnarray} and define the norms
\begin{eqnarray}
\nonumber&&  \norm{v}_V:=\norm{\nabla v}_{L^2}\,,\label{norm-V} \qquad \norm{q}_{P}=\norm{q}_{L_2}\,,\quad\qquad
\norm{X}_{Q}:=\norm{X}_{\mbox{\scriptsize H}(\mbox{\scriptsize
Curl};\Omega)},\label{norm-Q}\\
&&\norm{z}^2_{Z}:=\norm{v}^2_{V} +\norm{q}^2_{L^2}+\norm{X}^2_{Q}
+\norm{\beta}^2_{L^2}\quad\mbox{ for }
z=(v,q,X,\beta)\in \SFZ\,. \label{norm-Z}
\end{eqnarray}

Let us show that the bilinear form $\ba$ is coercive on $\SFW$.
 Let therefore $z=(v,q,X,\beta)\in \SFW$.
\begin{eqnarray}
\nonumber\ba(z,z)& \geq&  m_0\norm{\sym\nabla v-q}^2_{L^2}\mbox{
(from (\ref{ellipticityC}))} +\, \mu\,H_\chi\,\norm{\sym X-q}^2_{L^2}\\
\nonumber&&\qquad\qquad\hskip4.1truecm+\,\mu\,
L_c^2\,\norm{\Curl X}_{L^2}^2 +\mu\, k_2\,\norm{\beta}_{L^2}^2 \\
\nonumber&=&m_0\bigl[\norm{\sym\nabla v}^2_{L^2} +\norm{q}^2_{L^2}
-2\la\sym\nabla v, q\ra\bigr]+\mu\,H_\chi\,\bigl[\norm{\sym X}^2_{L^2} +\norm{q}^2_{L^2}\\
\nonumber&&\hskip5truecm-2\la\sym X, q\ra\bigr]  +\,\mu L_c^2\,\norm{\Curl X}_{L^2}^2
+\mu\,k_2\,\norm{\beta}_{L^2}^2\\
\nonumber&\geq &m_0\left[\norm{\sym\nabla v}^2_{L^2} +\norm{q}^2_{L^2} -\theta\norm{\sym\nabla
v}_{L^2}^2-\frac1\theta\norm{q}_{L^2}^2\right]\mbox{ (Young's inequality)}\\
\nonumber&& \hskip.5truecm+\, \mu\,H_\chi\,\left[\norm{\sym X}^2_{L^2} +\norm{
q}^2_{L^2} -\theta\norm{\sym X}_{L^2}^2-\frac1\theta\norm{q}_{L^2}^2\right]\mbox{ (Young's inequality)}\\
\nonumber&& \hskip.5truecm+\,\mu \,L_c^2\norm{\Curl X}_{L^2}^2+\frac12\,\mu\, k_2\,\norm{\beta}_{L^2}^2  +\,\frac12\,\mu\, k_2\,\norm{q}_{L^2}^2\mbox{ \,\,(using $\norm{q}\leq\beta$)} \\
\nonumber &=& m_0(1-\theta)\norm{\sym\nabla
v}^2_{L^2}+\left[(m_0+\mu\,H_\chi)\Bigl(1-\frac1\theta\Bigr)+\frac12\,\mu\,
k_2\right]\norm{q}_{L^2}^2\\
&&\hskip.5truecm+\,\,\mu\,H_\chi\,(1-\theta)\norm{\sym X}_{L^2}^2 +\,\mu \,L_c^2\,\norm{\Curl X}_{L^2}^2 +
\frac12\,\mu\,k_2\,\norm{\beta}_{L^2}^2\,.\label{coercive-a}\end{eqnarray}

So, since the hardening constant $k_2>0$, it is possible to choose  $\theta$ such that $$\displaystyle\frac{m_0+\mu\,H_\chi}{m_0+\mu\,H_\chi+\frac
12\,\mu\,k_2}< \theta<1,$$ and using Korn's first inequality (see e.g. \cite{NEFFKORN2002}) and the Korn-type inequality for incompatible tensor fields established  in \cite{NPW2011-1, NPW2012-1, NPW2012-2,NPW2014}, namely
\begin{equation}\label{korn-incomp}
\norm{X}^2_{L^2}\leq C(\norm{\sym X}^2_{L^2}+\norm{\Curl X}^2_{L^2})\qquad\forall X\in\mbox{H}_0(\mbox{Curl};\,\Omega,\,\Gamma_{\mbox{\scriptsize D}},\BBR^{3\times 3})\,,\end{equation}
there exists some constant
$C(m_0,\mu,H_\chi,k_2,L_c,\Omega)>0$ such that
\begin{equation}\label{coec-iso}a(z,z)\geq C\left[\norm{v}_V^2+\norm{q}^2_{L^2}+\norm{X}^2_{\mbox{\scriptsize H}(\mbox{\scriptsize
Curl};\Omega)}+\norm{\beta}_{L^2}^2
\right]=C\norm{z}^2_{Z}\quad\forall z=(v,q,X,\beta)\in \SFW\,.\end{equation} This shows existence for our microcurl model 
in polycrystalline gradient plasticity with isotropic hardening.

\begin{remark}\label{uniqueness-iso}{\rm Arguing as  in Section \ref{uniqueness}, we get for any two solutions $(u_i,\varepsilon_{p_{\,i}}, \PL_{p_{\,i}},\eta_{p_{\,i}})$ with $i=1,\,2$ of  (\ref{weak-form-iso-poly})  that 
$$u_1=u_2,\qquad\varepsilon_{p_{\,1}}=\varepsilon_{p_{\,2}},\qquad\eta_{p_{\,1}}=\eta_{p_{\,2}},\qquad\sym\PL_{p_{\,1}}=\sym\PL_{p_{\,2}},\qquad\Curl\PL_{p_{\,1}}=\Curl_{p_{\,2}}$$
Now, using the Korn-type inequality for incompatible tensor fields established  in \cite{NPW2011-1, NPW2012-1, NPW2012-2,NPW2014} and applied to $\PL_{p_{\,1}}-\PL_{p_{\,2}}$, namely
\begin{equation}\label{korn-incomp}
\norm{\PL_{p_{\,1}}-\PL_{p_{\,2}}}^2_{L^2}\leq C(\norm{\sym (\PL_{p_{\,1}}-\PL_{p_{\,2}})}^2_{L^2}+\norm{\Curl (\PL_{p_{\,1}}-\PL_{p_{\,2}})}^2_{L^2})\,,\end{equation} we also get that $\PL_{p_{\,1}}=\PL_{p_{\,2}}$ and this show the uniqueness of the weak/strong solution.
}\end{remark}

\subsection{The model with linear kinematical hardening}\label{kinhard}
Here we consider the model where the isotropic hardening has been replaced with linear kinematical hardening.  Here the free-energy  is given by
\begin{eqnarray}\label{freeng-poly-kin}
\Psi(\nabla u,\bvarepsilon_p,\PL_p,\Curl\PL_p):
&=&\underbrace{\Psi^{\mbox{\scriptsize lin}}_e(\bvarepsilon_e)}_{\mbox{\small elastic energy}}\qquad
+\qquad\underbrace{\Psi^{\mbox{\scriptsize lin}}_{\mbox{\scriptsize
 micro}}(\bvarepsilon_p,\PL_p)}_{\mbox{\small micro
energy}}\\
&& \nonumber+\,\,\underbrace{\Psi^{\mbox{\scriptsize lin}}_{\mbox{\scriptsize
curl}}(\Curl \PL_p)}_{\mbox{\small defect-like
energy (GND)}}\qquad+\hskip.4truecm\,\,\underbrace{\,\,\Psi^{\mbox{\scriptsize
lin}}_{\mbox{\scriptsize
 kin}}(\bvarepsilon_p)}_{\begin{array}{c}
\mbox{\small hardening energy (SSD)}\end{array}}\,,
 \end{eqnarray} where
 \begin{eqnarray}\label{freeng-expr-poly-kin} \nonumber\Psi^{\mbox{\scriptsize
lin}}_e(\bvarepsilon_e) &:=&\frac12\,\la\bvarepsilon_e,\C_{\mbox{\scriptsize iso}} \bvarepsilon_e\ra=\frac12\,\la\sym\nabla u-\bvarepsilon_p,\C_{\mbox{\scriptsize iso}}(\sym\nabla u-\bvarepsilon_p)\ra\,, \\
\Psi^{\mbox{\scriptsize lin}}_{\mbox{\scriptsize
 micro}}(\bvarepsilon_p,\PL_p)&:=&\frac12\,\mu\,H_\chi\,\norm{\bvarepsilon_p-\sym\PL_p}^2\,, \\ 
\nonumber \Psi^{\mbox{\scriptsize lin}}_{\mbox{\scriptsize
curl}}(\Curl \PL_p)&:=&\frac12\,\mu\, L_c^2\norm{\Curl \PL_p}^2\,,\qquad  \Psi^{\mbox{\scriptsize
lin}}_{\mbox{\scriptsize
 kin}}(\bvarepsilon_p)\,:=\,\frac12\mu\,k_1\,\norm{\bvarepsilon_p}^2
 \,.\end{eqnarray}

The strong formulation of the model is presented in Table \ref{table:micro-kinhard} while the weak formulation reads as

\begin{table}[ht!]{\footnotesize\begin{center}
\begin{tabular}{|ll|}\hline &\qquad\\
 {\em Additive split of strain:}&$\nabla u =e +p$,\quad $\bvarepsilon_e=\mbox{sym}\,e$,\quad $\bvarepsilon_p=\sym p$\\ 
{\em Equilibrium:} & $\mbox{Div}\,\sigma +f=0$ with
$\sigma=\C_{\mbox{\tiny iso}}\,\bvarepsilon_e=\C_{\mbox{\tiny iso}}(\sym\nabla u-\bvarepsilon_p)$\\
{\em Microbalance:} & $\mu\,L_c^2\,\Curl\Curl\PL_p=-\mu\,H_\chi\,(\sym\PL_p-\bvarepsilon_p)$,\\&\\
 {\em Free energy:} &
$\frac12\,\langle\C_{\mbox{\tiny iso}}\bvarepsilon_e,\bvarepsilon_e\rangle+\frac12\mu\,H_\chi\,\norm{\bvarepsilon_p-\sym\PL_p}^2$\\ &\qquad\qquad $+\frac12\,\mu\,
L^2_c\,\norm{\Curl \PL_p}^2+\frac12\,\mu\, k_1\,\norm{\sym p}^2$\\&\\
{\em Yield condition:} &
$\phi(\Sigma_E):=\norm{\dev\Sigma_E}-\yieldzero\leq0$\\
 {\em where } & $\Sigma_E:=\sigma+\Sigma^{\mbox{\tiny lin}}_{\mbox{\tiny micro}} 
+\Sigma^{\mbox{\tiny lin}}_{\mbox{\tiny kin}}$\\
{\em with }& $\Sigma^{\mbox{\tiny lin}}_{\mbox{\tiny micro}} =\mu\,H_\chi\,(\sym\PL_p-\bvarepsilon_p)=-\mu\,L_c^2\,\Curl\Curl\PL_p$\\
 & $\Sigma^{\mbox{\tiny lin}}_{\mbox{\tiny kin}}=-\mu\,k_1\bvarepsilon_p$
 \\&\\{\em Dissipation inequality:} &
 $\dsize\int_\Omega\langle\Sigma_E,\dot{\bvarepsilon}_p\ra\,dx\geq0$\\
 {\em Dissipation function:} &$\mathcal{D}_{\mbox{\tiny kin}}(q):=\yieldzero \norm{q}$\\
 {\em Flow law in primal form:} &
 $\Sigma_E\in\partial \mathcal{D}_{\mbox{\tiny kin}}(\dot{\bvarepsilon}_p)$\\&\\
{\em Flow law in dual form:}
&$\dot{\bvarepsilon}_p=\lambda\,\dsize\frac{\dev\Sigma_E}{\norm{\dev\Sigma_E}}$\\&\\
{\em KKT conditions:} &$\lambda\geq0$, \quad $\phi(\Sigma_E)\leq0$,
\quad $\lambda\,\phi(\Sigma_E)=0$\\
 {\em Boundary conditions for $\PL_p$:} & $\PL_p\times{n}=0$ on
 $\Gamma_{\mbox{\scriptsize D}}$,\,\, $(\Curl \PL_p)\times{n}=0$ on $\partial\Omega\setminus\Gamma_{\mbox{\scriptsize D}}$\\
 {\em Function space for $\PL_p$:} & $\PL_p(t,\cdot)\in \mbox{H}(\mbox{Curl};\;\Omega,\,\BBR^{3\times 3})$\\
 \hline
\end{tabular}\caption{\footnotesize The microcurl model in polycrystalline gradient plasticity with kinematical hardening.  Also in this model, the boundary condition on $\PL_p$ necessitates at least $\PL_p\in
\mbox{H}(\mbox{Curl};\,\Omega,\,\BBR^{3\times 3})$. Unlike the model with isotropic hardening for which uniqueness is obtained through the strong formulation, here we have uniqueness of the weak solution straight from the fomulation as a varational inequality.}\label{table:micro-kinhard}\end{center}}\end{table} 
\begin{eqnarray}\label{weak-form-iso}
\nonumber&&\hskip-.8truecm\int_\Omega\Bigl[\la\C_{\mbox{\tiny iso}}(\sym\nabla
u-\bvarepsilon_p),(\sym\nabla v-q)-(\sym\nabla\dot{u}-\dot{\bvarepsilon}_p)\ra+ \mu\,L_c^2\la\Curl\PL_p,\Curl\QL-\Curl\dot{\PL}_p\ra\\
\nonumber&&\qquad+\,\,\mu\,H_\chi\,\la\sym\PL_p-\bvarepsilon_p, (\sym\QL-q)-(\sym\dot{\PL}_p-\dot{\bvarepsilon}_p)\ra+\mu\,k_1\,\la\bvarepsilon_p\,q-\dot{\bvarepsilon}_p\ra\Bigr]dx\\
&&\qquad  +
\int_\Omega\mathcal{D}_{\mbox{\scriptsize kin}}(q)\,dx
-\int_\Omega\mathcal{D}_{\mbox{\scriptsize kin}}(\dot{\bvarepsilon}_p)\,dx\,\geq\, \int_\Omega
f\,(v-\dot{u})\,dx\,.\end{eqnarray}
That is, setting $\SFZ:=\SFV\times \SFP\times\SFQ$ with $\SFV$ and $\SFQ$ defined in (\ref{space-v})-(\ref{space-p}), $\SFP=L^2(\Omega,\mbox{Sym}(3)\cap\sL(3))$ and their norms in (\ref{norm-Z}), we get the problem of the form: Find $w=(u,\bvarepsilon_p,\PL_p)\in  \SFH^1(0,T;Z)$
such that $w(0)=0$ and 
\begin{equation}\label{wf-kin}
\ba(\dot{w},z-w)+j(z)-j(\dot{w})\geq \langle
\ell,z-\dot{w}\rangle\mbox{ for every } z\in \SFZ\mbox{ and for a.e.
}t\in[0,T]\,,\end{equation}
where
\begin{eqnarray}
\ba(w,z)&=&\int_\Omega\Bigl[\la\C_{\mbox{\tiny iso}}(\sym\nabla
u-\bvarepsilon_p),\sym\nabla v-q\ra+ \mu\,L_c^2\la\Curl\PL_p,\Curl X\ra\\
\nonumber&&\hskip2truecm+\,\,\mu\,H_\chi\,\la\sym\PL_p-\bvarepsilon_p, \sym X-q\ra+\mu\,k_1\,\la\bvarepsilon_p,q\ra\,\Bigr]dx\,,
\label{bilin-iso-spin}
\\\nonumber\\
j(z)&=&\int_\Omega\mathcal{D}_{\mbox{\tiny kin}}(q)\,dx\,,\label{functional-kin}\\
\langle \ell,z\rangle&=&\int_\Omega
f\,v\,dx\,,\label{lin-form-kin}\end{eqnarray} for
$w=(u,\bvarepsilon_p,\PL_p)$ and $z=(v,q,X)$ in
$\SFZ$.
\vskip.2truecm\noindent The existence and uniqueness result for the problem in  (\ref{wf-kin})-(\ref{lin-form-kin}) is obtained from  \cite[Theorem 6.15]{Han-ReddyBook} as the bilinear form $\ba$ is $\SFZ$-coercive (arguing as in 
 (\ref{coercive-a}).
 
 \section{The relaxed linear micromorphic continuum}\label{micromorphic}
 The relaxed micromorphic model is a very special subclass of the micromorphic model approach in which the extra dependence on gradients of the micro-distortion appears only through the Curl-operator. In the static and isotropic cases, the purely elastic model consists of a two-field minimization problem for the displacement $u:\Omega\subset\mathbb{R}^3\to\mathbb{R}^3$ and the non-symmetric micro-distortion tensor $\PL_p:\Omega\subset\mathbb{R}^3\to\mathbb{R}^{3\times 3}$ so that for
 \begin{eqnarray}\label{micromorphic1}
\nonumber  \mathcal{E}(u,\PL_p)&:=&\int_\Omega\Bigl[\mu_e\norm{\sym(\nabla u-\PL_p)}^2 +\mu_c\norm{\skew(\nabla u-\PL_p)}^2+\frac{\lambda_e}2(\tr(\nabla u-\PL_p))^2\\
&& \qquad+\,\mu_{\mbox{\tiny micro}}\norm{
 \sym\PL_p}^2+\frac{\lambda_{\mbox{\tiny micro}}}2(\tr(\PL_p))^2 +\mu_e\,\frac{L_c^2}2\,\norm{\Curl\PL_p}^2\Bigr]\,dx\,,\end{eqnarray}
 \begin{equation}\label{Eq1-micromorphic}
 \mathcal{E}(u,\PL_p)\to\min\mbox{ w.r.t }(u,\,\PL_p)
 \end{equation} subject to displacement boundary conditions $u|_{\Gamma_{\mbox{\scriptsize D}}}=0$ and the tangential boundary conditions $\PL_p\times n|_{\Gamma_{\mbox{\scriptsize D}}}=0$ (equivalent to $\PL_p\cdot\tau|_{\Gamma_{\mbox{\scriptsize D}}}=0$ for all vectors $\tau$ tangent to $\Gamma_{\mbox{\scriptsize D}}$).
 Here, $\mu_e$ and $\lambda_e$ with 
 \begin{equation}\label{Lame1}
 \mu_e>0\quad\mbox{ and }\quad 2\mu_e+3\lambda_e>0\,,\end{equation} are new elastic material constants which are {\bf not} the Lam\'e constants of linear elasticity. Well-posedness results in statics and dynamics have been obtained in \cite{NGMPR2014, NGMP2014}, making crucial use of a recently established Korn's inequality for incompatible tensor fields \cite{NPW2011-1, NPW2012-1, NPW2012-2,NPW2014}.  The parameter $\mu_c\geq0$ is called the Cosserat couple modulus and may be set to zero in this model.
 
 Regarding the relation to the polycrystalline  microcurl model (\ref{freeng-poly-kin})-(\ref{freeng-expr-poly-kin}), we see that in (\ref{micromorphic1}) the minimization variable $\PL_p$ is elastically coupled to the displacement gradient $\nabla u$ instead of being (penalty)-coupled to the plastic distortion $p$ in the microcurl model  (\ref{freeng-poly-kin})-(\ref{freeng-expr-poly-kin}).
 
  In the single crystal microcurl model, the equation for the micro-distortion can be obtained from the one-field minimization problem
 \begin{equation}\label{Eq2-micromorphic}
 \int_\Omega\Bigl[\frac12\,\mu\,H_\chi\,\norm{p-\PL_p}^2+\frac12\,\mu\,L_c^2\,\norm{\Curl\PL_p}^2\Bigr]\,dx\,\,\to\,\,\min\,\,\PL_p
 \end{equation} at given plastic distortion $p$.
 
 Now, if we let $\mu_e,\,\mu_c,\,\lambda_e\to\infty$ ($\PL_p\to\nabla u$) then the static model turns indeed into a linear elastic model 
 \begin{equation}\label{lin-elast-model}
 \int_\Omega\Bigl[\mu_\infty\,\norm{\sym\nabla u}^2+\frac{\lambda_\infty}2\,(\tr(\nabla u))^2\Bigr]\,dx\,\,\to\,\,\min\,\,u\,,\end{equation} where $\lambda_\infty,\,\,\mu_\infty$ can be determined analytically \cite{BDAGMP2016}.
 
 The formulation (\ref{micromorphic1}) in the dynamic case has a number of distinguishing features. As it turns out, the so-callled metamaterials with band-gaps at certain frequency ranges can be qualitatively and quantatively described. For this, a nonzero Cosserat couple modulus $\mu_c>0$ is mandatory. Materials that do not show band-gaps must be modelled with $\mu_c=0$.   Note that the formulation (\ref{micromorphic1}) contains as the special case $\mu_{\mbox{\tiny micro}},\,\,\lambda_{\mbox{\tiny micro}}\to\infty$ the well-known infinitesimal Cosserat model in which the additional field $\PL_p$ is restricted to be skew-symmetric (i.e., $\PL_p$ is set as $A\in\so(3)$) and the elastic minimization problem reads 
 \begin{eqnarray}\label{lin-cosserat1}
 &&\int_\Omega\Bigl[\mu\,\norm{\sym\nabla u}^2+\frac\lambda2\,(\tr(\nabla u))^2+\mu_c\,\norm{\skew(\nabla u-A)}^2+ \mu\,\frac{L_c^2}2\,\norm{\Curl A}^2\Bigr]\,dx\\
 \nonumber&&\qquad\,\to\,\min\,\,(u,A)\,, 
 \end{eqnarray}
 see e.g. \cite{BDAGMP2016}. The latter formulation has been coupled to perfect plasticity in an endevour to regularize ill-posedness of perfect plasticity, see e.g. \cite{neffchel2007, NCMW2007}

\section{Conclusion}
Examples of finite element computations based on the microcurl single crystal models can be found in \cite{CORFORBUS}
where polycrystalline microstructures are discetized
in order to account for grain size effects on the local stress and lattice curvatures fields inside the grains and 
on the overall Hall-Petch effect.
Orowan-type size effects were addressed for laminate microstructures in \cite{WULFORBOH2015}. It remains to implement the polycrystalline
formulation proposed in the present work and to compare its response to that of polycrystalline aggregates using the single crystal model.
In that way the new material parameters could be identified from this multiscale analysis. This would also help to decide between the two possible penalty couplings, namely 
$$\underbrace{\frac12\,\mu\,H_\chi\norm{p-\PL_p}^2}_{\mbox{direct coupling}}\qquad\qquad\mbox{versus}\qquad\qquad\underbrace{\frac12\,\mu\,H_\chi\norm{\sym(p-\PL_p)}^2}_{\begin{array}{c}
\mbox{symmetric coupling}\\\mbox{treated in this work}\end{array}}\,.$$
 Mathematically, both formulations are well-posed, provided sufficient hardening is present. The direct coupling has the advantage of a clear penalty interpretation while the symmetric coupling does not see the plastic spin altogether, which may be advantageous from a modelling and implementational point of view.
 
The present mathematical analysis was performed within the infinitesimal framework. The reader is referred to \cite{acgferingenvol11}
for a finite deformation formulation of the microcurl single crystal model that can be used for further applications involving significant 
lattice rotations and strains.

\end{document}